\documentclass[10pt,twocolumn]{IEEEtran}
\usepackage[utf8]{inputenc}
\usepackage{mathrsfs}
\usepackage{amsmath}
\usepackage{amsthm}
\usepackage{amssymb}
\usepackage{cancel}
\usepackage[unicode=true,
 bookmarks=false,
 breaklinks=false,pdfborder={0 0 0},pdfborderstyle={},backref=false,colorlinks=false]
 {hyperref}

\makeatletter

\providecommand{\tabularnewline}{\\}

\theoremstyle{plain}
\newtheorem{thm}{\protect\theoremname}
\theoremstyle{definition}
\newtheorem{example}[thm]{\protect\examplename}
\theoremstyle{plain}
\newtheorem{prop}[thm]{\protect\propositionname}
\theoremstyle{plain}
\newtheorem{lem}[thm]{\protect\lemmaname}

\IEEEoverridecommandlockouts
\usepackage{comment}
\usepackage{tikz}
\usepackage{graphicx}
\usepackage{amsfonts}
\usepackage{mathrsfs}
\usepackage{fix2col}
\usepackage{latexsym}
\usepackage[dvips]{epsfig}
\usepackage{cite}
\usepackage{graphics}
\usepackage{epstopdf}
\usepackage{enumerate}
\usepackage{cite}
\usepackage{epstopdf}
\usepackage{arydshln}

\usepackage{color}

\newtheorem{assume}{Assumption}

\makeatother

\providecommand{\examplename}{Example}
\providecommand{\lemmaname}{Lemma}
\providecommand{\propositionname}{Proposition}
\providecommand{\theoremname}{Theorem}

\begin{document}

\title{\Large \bf Conic Optimization Theory: Convexification Techniques
and\\
 Numerical Algorithms \thanks{Email: ryz@berkeley.edu, cedric.josz@gmail.com, sojoudi@berkeley.edu}}

\author{Richard Y. Zhang\textsuperscript{{*}}, C\'edric Josz\textsuperscript{{*}},
and Somayeh Sojoudi\\
 \thanks{\protect\textsuperscript{{*}} The first two authors contributed equally
to this work.} \thanks{Richard Y. Zhang is with the Department of Industrial Engineering
and Operations Research, University of California, Berkeley. Cedric
Josz is with the Department of Electrical Engineering and Computer
Sciences, University of California, Berkeley. Somayeh Sojoudi is with
the Departments of Electrical Engineering and Computer Sciences and
Mechanical Engineering, University of California, Berkeley. This work
was supported by the ONR grant N00014-17-1-2933, DARPA grant D16AP00002,
and AFOSR grant FA9550-17-1-0163.}}
\maketitle
\begin{abstract}
Optimization is at the core of control theory and appears in several
areas of this field, such as optimal control, distributed control,
system identification, robust control, state estimation, model predictive
control and dynamic programming. The recent advances in various topics
of modern optimization have also been revamping the area of machine
learning. Motivated by the crucial role of optimization theory in
the design, analysis, control and operation of real-world systems,
this tutorial paper offers a detailed overview of some major advances
in this area, namely conic optimization and its emerging applications.
First, we discuss the importance of conic optimization in different
areas. Then, we explain seminal results on the design of hierarchies
of convex relaxations for a wide range of nonconvex problems. Finally,
we study different numerical algorithms for large-scale conic optimization
problems.
\end{abstract}

\section{Introduction}

Optimization is an important tool in the design, analysis, control
and operation of real-world systems. In its purest form, optimization
is the mathematical problem of minimizing (or maximizing) an \emph{objective
function} by selecting the best choice of \emph{decision variables},
possibly subject to \emph{constraints} on their specific values. In
the wider engineering context, optimization also encompasses the process
of identifying the most suitable objective, variables, and constraints.
The goal is to select a mathematical model that gives useful insight
to the practical problem at hand, and to design a robust and scalable
algorithm that finds a provably optimal (or near-optimal) solution
in a reasonable amount of time. 

The theory of \emph{convex} optimization had a profound influence
on the development of modern control theory, giving rise to the ideas
of robust and optimal control~\cite{boyd1994linear,zhou1996robust,dullerud2000course},
distributed control~\cite{boyd2011distributed}, system identification~\cite{ljung1998system},
model predictive control~\cite{camacho2013model} and dynamic programming~\cite{bertsekas1995dynamic}.
The mathematical study of convex optimization dates by more than a
century, but its usefulness for practical applications only came to
light during the 1990s, when it was discovered that many important
engineering problems are actually convex (or can be well-approximated
as being so)~\cite{boyd2004convex}. Convexity is crucial in this
regard, because it allows the corresponding optimization problem to
be solved—very reliably and efficiently—to global optimality, using
only local search algorithms. Moreover, this global optimality can
be certified by solving a dual convex optimization problem. 

With the recent advances in computing and sensor technologies, convex
optimization has also become the backbone of data science and machine
learning, where it supplies us the techniques to extract useful information
from data~\cite{bach2004multiple,buhlmann2011statistics,Simon13}.
This tutorial paper is, in part, inspired by the crucial role of optimization
theory in both the long-standing area of control systems and the newer
area of machine learning, as well as its multi-billion applications
in large-scale real-world systems such as power grids. 

Nevertheless, most interesting optimization problems are nonconvex:
structured analysis and synthesis and output feedback control in control
theory~\cite{zhou1996robust,parrilo2000structured}; deep learning,
Bayesian inference, and clustering in machine learning~\cite{bach2004multiple,lanckriet2004learning,weinberger2006unsupervised,boyd2011distributed};
and integer and mixed integer programs in operations research~\cite{conforti2014,kocuk2017new,fattahi2017conic}.
Nonconvex problems are difficult to solve, both in theory and in practice.
While candidate solutions are easy to find, locally optimal solutions
are not necessarily globally optimal. Even if a candidate solution
is suspected of being globally optimal, there is often no effective
way of validating the hypothesis. 

Convex optimization can rigorously solve nonconvex problems to global
optimality, using techniques collectively known as \emph{convexification}.
The essential idea is to relax a nonconvex problem into a hierarchy
of convex problems that monotonously converge towards the global solution~\cite{parrilo-2000b,lasserre-2001,laurent-2003,lasserre-2010}.
These resulting convex problems come in a standard form known as \emph{conic
optimization}, that generalizes semidefinite programs (SDP) and linear
matrix inequalities (LMI) from control theory, as well as quadratically-constrained
quadratic programs (QCQPs) from statistics, and linear programs (LPs)
from operations research. Conic optimization can be solved with reliability
using interior-point methods~\cite{alizadeh1995interior,nesterov1997self},
and also with great efficiency using large-scale numerical algorithms
designed to exploit problem structure~\cite{fukuda2001exploiting,toh2002solving,burer2003nonlinear,kovcvara2003pennon,journee2010low,zhao2010newton,zhang2017modified,nesterov2007smoothing,vandenberghe2015chordal}.

The remainder of this paper is organized as follows. Case studies
are provided in Section~\ref{sec:2} to show the importance of conic
optimization in emerging applications (other than well-known problems
in control theory). Different convexification techniques for polynomial
optimization are studied in Section~\ref{sec:3}, followed by a detailed
investigation of numerical algorithms for conic optimization in Section~\ref{sec:4}.
Concluding remarks are drawn in Section~\ref{sec:5}.

\textbf{Notations:} The symbols $\mathbb{R}$ and $\mathbb{S}^{n}$
denote the sets of real numbers and $n\times n$ real symmetric matrices,
respectively. The symbols $\mathbb{R}_{+}$ and $\mathbb{S}_{+}^{n}$
denote the sets of nonnegative real numbers and $n\times n$ symmetric
and positive-semidefinite matrices, respectively. $\text{rank}\{\cdot\}$
and $\text{trace}\{\cdot\}$ denote the rank and trace of a matrix.
The notation $X\succeq0$ means that $X\in\mathbb{S}_{+}^{n}$. $X^{\text{opt}}$
shows a global solution of a given conic optimization problem with
the variable $X$.

The vectorization of a matrix is the column-stacking operation 
\[
\mathrm{vec}\,X=[X_{1,1},\ldots,X_{n,1},X_{1,2},\ldots,X_{n,2},\ldots,X_{n,n}]^{T},
\]
and the Kronecker product 
\[
A\otimes B=\begin{bmatrix}A_{1,1}B & \cdots & A_{1,n}B\\
\vdots & \ddots & \vdots\\
A_{n,1}B & \cdots & A_{n,n}B
\end{bmatrix},
\]
is defined to satisfy the Kronecker identity 
\[
\mathrm{vec}\,(AXB)=(B^{T}\otimes A)\mathrm{vec}\,X.
\]

\section{Emerging Applications of Conic Optimization}

\label{sec:2}

Although conic optimization has appeared in many subareas of control
theory since early 1990s, it has found new applications in many other
problems in the last decade. Some of these applications will be discussed
below.

\subsection{Machine Learning}

In machine learning, kernel methods are used to study relations, such
as principle components, clusters, rankings, and correlations, in
large datasets~\cite{bach2004multiple}. They are used for pattern
recognition and analysis. Support vector machines (SVMs) are kernel-based
supervised learning techniques. SVMs are one of the most popular classifiers
that are currently used. Kernel matrix is a symmetric and positive
semidefinite matrix that plays a key role in kernel-based learning
problems. Specifying this matrix is a classical model selection problem
in machine learning. A great deal of effort has been made over the
past two decades to elucidate the role of semidefinite programing
as an efficient convex optimization framework for machine learning,
kernel-machines, SVMs, and learning kernel matrices from data in particular~\cite{lanckriet2004learning,bach2004multiple,kato2008multi,graepel2002kernel,weinberger2006unsupervised,graepel2004invariant}.

In many applications, such as social networks, neuroscience, and financial
markets, there exist a massive amount of multivariate timestamped
observations. Such data can often be modeled as a network of interacting
components, in which every component in the network is a node associated
with some time series data. The goal is to infer the relationships
between the network entities using observational data. Learning and
inference are important topics in machine learning. Learning a hidden
structure in data is usually formulated as an optimization problem
augmented with sparsity-promoting techniques~\cite{Coleman90,Bach11,Benson}.
These techniques have become essential to the tractability of big-data
analyses in many applications, such as data mining~\cite{Garcke01,Muth05,Wu14},
pattern recognition~\cite{Wright10,Qiao10}, human brain functional
connectivity~\cite{Sojoudi14}, and compressive sensing~\cite{Candes07,Simon13}.
Similar approaches have been used to arrive at a parsimonious estimation
of high-dimensional data. However, most of the existing statistical
learning techniques in data analytics need a large amount of data
(compared to the number of parameters), which limit their applications
in practice~\cite{buhlmann2011statistics,fan2010selective}. To address
this issue, a special attention has been paid to the augmentation
of these learning problems with sparsity-inducing penalty functions
to obtain sparse and easy-to-analyze solutions.

Graphical lasso (GL) is a popular method for estimating the inverse
covariance matrix~\cite{friedman2008sparse,banerjee2008model,yuan2007model}.
GL is an optimization problem that shrinks the elements of the inverse
covariance matrix towards zero compared to the maximum likelihood
estimates, using an $l_{1}$ regularization. There is a large body
of literature suggesting that the solution of GL is a good estimate
for the unknown graphical model, under a suitable choice of the regularization
parameter \cite{friedman2008sparse,banerjee2008model,yuan2007model,liu2010stability,kramer2009regularized,danaher2014joint}.
To explain the GL problem, consider a random vector $x=(x_{1},x_{2},...,x_{d})$
with a multivariate normal distribution. Let $\Sigma_{*}\in\mathbb{S}^{d}$
denote the correlation matrix associated with the vector $x$. The
inverse of the correlation matrix can be used to determine the conditional
independence between the random variables $x_{1},x_{2},...,x_{d}$.
In particular, if the $(i,j)^{\text{th}}$ entry of $\Sigma_{*}^{-1}$
is zero for two disparate indices $i$ and $j$, then $x_{i}$ and
$x_{j}$ are conditionally independent given the rest of the variables.
The graph $\text{supp}\big(\Sigma_{*}^{-1}\big)$ (i.e., the sparsity
graph of $\Sigma_{*}^{-1}$) represents a graphical model capturing
the conditional independence between the elements of $\bold x$.

Assume that $\Sigma_{*}$ is nonsingular and that $\text{supp}\big(\Sigma_{*}^{-1}\big)$
is a sparse graph. Finding this graph is cumbersome in practice because
the exact correlation matrix $\Sigma_{*}$ is rarely known. More precisely,
$\text{supp}\big(\Sigma_{*}^{-1}\big)$ should be constructed from
a given sample correlation matrix (constructed from $n$ samples),
as opposed to $\Sigma_{*}$. Let $\Sigma$ denote an arbitrary $d\times d$
positive-semidefinite matrix, which is provided as an estimate of
$\Sigma_{*}$. Consider the convex optimization problem: 
\begin{equation}
\begin{aligned}\min_{S\in\mathbb{S}^{d}}\quad & -\log\det(S)+\mathrm{trace}(\Sigma S)\\
\text{s.t.}\quad & \quad S\succeq0
\end{aligned}
\label{eq1_o}
\end{equation}
It is easy to verify that the optimal solution of the above problem
is equal to $S^{\text{opt}}=\Sigma^{-1}$. However, there are two
issues with this solution. First, since the number of samples available
in many applications is modest compared to the dimension of $\Sigma$,
the matrix $\Sigma$ could be ill-conditioned or even singular. In
that case, the equation $S^{\text{opt}}=\Sigma^{-1}$ leads to large
or unbounded entries for the optimal solution of \eqref{eq1_o}. Second,
although $\Sigma_{*}^{-1}$ is assumed to be sparse, a small random
difference between $\Sigma_{*}$ and $\Sigma$ would make $S^{\text{opt}}$
highly dense. In order to address the aforementioned issues, consider
the problem 
\begin{equation}
\begin{aligned}\min_{S\in\mathbb{S}^{d}}\quad & -\log\det(S)+\mathrm{trace}(\Sigma S)+\lambda\|S\|_{*}\\
\text{s.t.}\quad & \quad S\succeq0
\end{aligned}
\label{eq1}
\end{equation}
where $\lambda\in\mathbb{R}_{+}$ is a regularization parameter and
$\|S\|_{*}$ denotes the sum of the absolute values of the off-diagonal
entries of $S$. This problem is referred to as \textit{Graphical
Lasso} (GL). Intuitively, the term $\|S\|_{*}$ in the objective function
serves as a surrogate for promoting sparsity among the off-diagonal
entries of $S$, while ensuring that the problem is well-defined even
with a singular input $\Sigma$.

There have been major interests in studying the properties of GL as
a conic optimization problem, in addition to the design of numerical
algorithms for this problem \cite{friedman2008sparse,banerjee2008model,Hsieh14}.
For example, \cite{Sojoudi16} and \cite{Somayeh14} have shown that
the conic optimization problem~\eqref{eq1} is highly related to
a simple thresholding technique. This result is leveraged in \cite{fattahi2017graphical}
to obtain an explicit formula that serves as an exact solution of
GL for acyclic graphs and as an approximate solution of GL for arbitrary
sparse graphs. Another line of work has been devoted to studying the
connectivity structure of the optimal solution of the GL problem.
In particular, \cite{Mazumdar12} and \cite{Witten11} have proved
that the connected components induced by thresholding the sample correlation
matrix and the connected components in the support graph of the optimal
solution of the GL problem lead to the same vertex partitioning.

\subsection{Optimization for Power Systems}

The real-time operation of an electric power network depends heavily
on several large-scale optimization problems solved from every few
minutes to every several months. State estimation, optimal power flow
(OPF), unit commitment, and network reconfiguration are some fundamental
optimization problems solved for transmission and distribution networks.
These different problems have all been built upon the power flow equations.
Regardless of their large-scale nature, it is a daunting challenge
to solve these problems efficiently. This is a consequence of the
nonlinearity/non-convexity created by two different sources: (i) discrete
variables such as the ratio of a tap-changing transformer, the on/off
status of a line switch, or the commitment parameter of a generator,
and (ii) the laws of physics. Issue (i) is more or less universal
and researchers in many fields of study have proposed various sophisticated
methods to handle integer variables. In contrast, Issue (ii) is pertinent
to power systems, and it demands new specialized techniques and approaches.
More precisely, complex power being a quadratic function of complex
bus voltages imposes quadratic constraints on OPF-based optimization
problems, and makes them NP-hard~\cite{Bienstock15}.

OPF is at the heart of Independent System Operator (ISO) power markets
and vertically integrated utility dispatch~\cite{FERC_1}. This problem
needs to be solved annually for system planning, daily for day-ahead
commitment markets, and every 5-15 minutes for real-time market balancing.
The existing solvers for OPF-based optimization either make potentially
very conservative approximations or deploy general-purpose local-search
algorithms. For example, a linearized version of OPF, named DC OPF,
is normally solved in practice, whose solution may not be physically
meaningful due to approximating the laws of physics. Although OPF
has been studied for 50 years, the algorithms deployed by ISOs suffer
from several issues, which may incur tens of billions of dollars annually~\cite{FERC_1}.

The power flow equations for a power network are quadratic in the
complex voltage vector. It can be verified that the constraints of
the above-mentioned OPF-based problems can often be cast as quadratic
constraints after introducing certain auxiliary parameters. This has
inspired many researchers to study the benefits of conic optimization
for power optimization problems. In particular, it has been shown
in a series of papers that a basic SDP relaxation of OPF is exact
and finds global minima for benchmarks examples~\cite{LL2012_1,LZT2012_1,madaniconvex2013,sojoudiconvexification2013}.
By leveraging the physics of power grids, it is also theoretically
proven that the SDP relaxation is always exact for every distribution
network and every transmission network containing a sufficient number
of transformers (under some technical assumptions)~\cite{SL2012_1,javadsojoudi14}.
\begin{table}
\caption{Performance of penalized SDP for OPF.}
\label{case-stu} 

\centering %
\begin{tabular}{|l|c|c|c|}
\hline 
\hspace{5mm} \textbf{{Test }} & \textbf{{Near-optimal }} & \textbf{{Global optimality}} & \textbf{{Run}}\tabularnewline
\hspace{4mm} \textbf{{ cases} } & \textbf{{cost}} & \textbf{{guarantee}} & \textbf{{time (s)}}\tabularnewline
\hline 
Polish 2383wp  & 1874322.65 & 99.316\%  & 529\tabularnewline
\hline 
Polish 2736sp  & 1308270.20 & 99.970\% & 701 \tabularnewline
\hline 
Polish 2737sop  & 777664.02 & 99.995\% & 675\tabularnewline
\hline 
Polish 2746wop  & 1208453.93 & 99.985\%  & 801\tabularnewline
\hline 
Polish 2746wp  & 1632384.87 & 99.962\%  & 699\tabularnewline
\hline 
Polish 3012wp  & 2608918.45 & 99.188\% & 814 \tabularnewline
\hline 
Polish 3120sp  & 2160800.42 & 99.073\%  & 910\tabularnewline
\hline 
\end{tabular} 
\end{table}

The papers~\cite{Madani_2014_Allerton_1} and~\cite{madaniconvex2013}
show that if the SDP relaxation is not exact due to the violation
of certain assumptions, a penalized SDP relaxation would work for
a carefully chosen penalty term, which leads to recovering a near-global
solution. This technique is tested on several real-world grids and
the outcome is partially reported in Table~\ref{case-stu}~\cite{OPF_Solver_2014}.
It can be observed that the SDP relaxation has found operating points
for the nationwide grid of Poland in different times of the year,
where the global optimality guarantee of each solution is at least
99\%, implying that the unknown global minima are at most $1\%$ away
from the obtained feasible solutions. In some cases, finding a suitable
penalization parameter can be challenging, but this can be remedied
by using~\cite[Algorithm 1]{laplacian_obj} based on the notion of
Laplacian matrix. Different techniques have been proposed in the literature
to obtain tighter relaxations of OPF~\cite{chen2016bound,coffrin2015,coffrin2015bis,ghaddar2016,kocuk2015}.
Several papers have developed SDP-based approximations or reformulations
for a wide range of power problems, such as state estimation~\cite{weng2015convexification,HICSS_Ramtin_2017},
unit commitment~\cite{fattahi2017conic}, and transmission switching~\cite{kocuk2017new,fattahi_OTS_2017}.
The recent findings in this area show the significant potential of
conic relaxation for power optimization problems.

\subsection{Matrix Completion}

Consider a symmetric matrix $\widehat{X}\in\mathbb{S}^{n}$, where
certain off-diagonal entries are missing. The low-rank positive-semidefinite
matrix completion problem is concerned with the design of the unknown
entries of this matrix in such a way that the matrix becomes positive
semidefinite with the lowest rank possible. This problem has applications
in signal processing and machine learning, where the goal is to learn
a structured dataset (or signal) from limited observations (or samples).
It is also related to the complexity reduction for semidefinite programs~\cite{fukuda2001exploiting,nakata2003exploiting,grone1984positive}.
Let the known entries of $\widehat{X}$ be represented by a graph
$\mathcal{G}=(\mathcal{V},\mathcal{E})$ with the vertex set $\mathcal{V}$
and the edge set $\mathcal{E}$, where each edge of the graph corresponds
to a known off-diagonal entry of $\widehat{X}$. The matrix completion
problem can be expressed as: \begin{subequations} \label{eq_nn_1}
\begin{align}
\min_{X\in\mathbb{S}^{n}}\quad\quad & \text{rank}\{X\}\\
\text{s.t.}\hspace{0.8cm} & X_{ij}=\widehat{X}_{ij},\qquad\quad\ \forall(i,j)\in\mathcal{E}\\
 & X_{kk}=\widehat{X}_{kk},\qquad\quad\forall k\in\mathcal{V}\\
 & X\succeq0
\end{align}
\end{subequations} The missing entries of $\widehat{X}$ can be obtained
from an optimal solution of the above problem. Since~\eqref{eq_nn_1}
is nonconvex, it is desirable to find a convex formulation/approximation
of this problem. To this end, consider a number $\bar{n}$ that is
greater than or equal to $n$. Let $\bar{\mathcal{G}}=(\bar{\mathcal{V}},\bar{\mathcal{E}})$
be a graph with $\bar{n}$ vertices such that $\mathcal{G}$ is a
subgraph of $\bar{\mathcal{G}}$. Consider the optimization problem
\begin{subequations} \label{eq_nn_2} 
\begin{align}
\min_{\bar{X}\in\mathbb{S}^{\bar{n}}}\quad\quad & \hspace{-3mm}\sum_{(i,j)\in\bar{\mathcal{E}}\backslash\mathcal{E}}t_{ij}\ {\bar{X}}_{ij}\\
\ \text{s.t.}\hspace{0.8cm} & \bar{X}_{ij}=\widehat{X}_{ij},\qquad\quad\ \forall(i,j)\in\mathcal{E}\\
 & \bar{X}_{kk}=\widehat{X}_{kk},\qquad\quad\forall k\in\mathcal{V}\\
 & \bar{X}_{kk}=1,\qquad\quad\hspace{0.5cm}\forall k\in\bar{\mathcal{V}}\backslash\mathcal{V}\\
 & \bar{X}\succeq0
\end{align}
\end{subequations} Define $\bar{X}^{\text{opt}}(n)$ as the $n\times n$
principal submatrix of an optimal solution of~\eqref{eq_nn_2}. It
is shown in \cite{Madani_SIAM_2017} that: 
\begin{itemize}
\item $\bar{X}^{\text{opt}}(n)$ is a positive semidefinite filling of $\widehat{X}$
whose rank is upper bounded by certain parameters of the graph $\bar{\mathcal{G}}$,
no mater what the coefficients $t_{ij}$'s are as long as they are
all nonzero. 
\item $\bar{X}^{\text{opt}}(n)$ is low rank if the graph $\mathcal{G}$
is sparse.
\item $\bar{X}^{\text{opt}}(n)$ is guaranteed to be a solution of~\eqref{eq_nn_1}
for certain types of graphs $\mathcal{G}$.
\end{itemize}
Since~\eqref{eq_nn_2} is a semidefinite program, it points to the
role of conic optimization in solving the matrix completion problem.

\subsection{Affine Rank Minimization Problem}

\label{subsect:affine}

Consider the problem \begin{subequations} \label{eq_nn_3} 
\begin{align}
\min_{Y\in\mathbb{R}^{m\times r}}\quad\quad & \text{rank}\{Y\}\\
\text{s.t.}\hspace{0.8cm} & \text{trace}\{{N}_{k}{Y}\}\leq a_{k}, & k=1,\dots,p
\end{align}
\end{subequations} where ${N}_{1},\ldots,{N}_{p}\in\mathbb{R}^{r\times m}$
are constant sparse matrices. This is a general affine rank minimization
problem without any positive semidefinite constraint. A special case
of the above problem is the regular matrix completion problem defined
as finding the missing the entries of a rectangular matrix with partially
known entries to minimize its rank \cite{johnson1990matrix,candes2009exact,recht2010guaranteed,keshavan2010matrix,candes2010power}.
A popular convexification method for~\eqref{eq_nn_3} is to replace
its objective function with the nuclear norm of $Y$~\cite{recht2010guaranteed}.
This is motivated by the fact that the nuclear norm function is the
convex envelop of $\text{rank}\{Y\}$ on the set $\{Y\in\mathbb{R}^{m\times r}\,\vert\,\|Y\|\leq1\}$~\cite{fazel2002matrix}.
Optimization \eqref{eq_nn_3} can be reformulated as a matrix optimization
problem whose matrix variable is symmetric and positive semidefinite.
To explain this reformulation, consider a new matrix variable $X$
defined as 
\begin{equation}
{X}\triangleq\left[\begin{array}{cc}
Z_{1} & Y\\
Y^{\mathrm{T}} & Z_{2}
\end{array}\right],
\end{equation}
where $Z_{1}$ and $Z_{2}$ are auxiliary matrices, and $Y$ acts
as a submatrix of $X$ corresponding to its first $m$ rows and last
$r$ columns. Now, consider the problem \begin{subequations} \label{eq_nn_4}
\begin{align}
\min_{X\in\mathbb{S}^{m+r}}\quad\quad & \text{rank}\{X\}\\
\text{s.t.}\hspace{0.9cm} & \text{trace}\{{M}_{k}{X}\}\leq a_{k}, & k=1,\dots,p\\
 & X\succeq0
\end{align}
\end{subequations} where 
\begin{equation}
{M}_{k}\triangleq\left[\begin{array}{cc}
0_{m\times m} & \frac{1}{2}{N}_{k}^{\mathrm{T}}\\
\frac{1}{2}{N}_{k} & 0_{r\times r}
\end{array}\right]
\end{equation}
For every feasible solution ${X}$ of the above problem, its associated
submatrix $Y$ is feasible for \eqref{eq_nn_3} and satisfies the
inequality 
\begin{equation}
\text{rank}\{Y\}\leq\text{rank}\{X\}
\end{equation}
The above inequality turns into an equality at optimality and, moreover,
the problems \eqref{eq_nn_3} and \eqref{eq_nn_4} are equivalent
\cite{fazel2002matrix,fazel2003log}. One may replace the rank objective
in \eqref{eq_nn_4} with the linear term $\text{trace}\{X\}$ based
on the nuclear norm technique, or use the more general idea delineated
in \eqref{eq_nn_2} to find an SDP approximation of the problem with
a guarantee on the rank of the solution.

\subsection{Conic Relaxation of Quadratic Optimization}

Consider the standard non-convex quadratically-constrained quadratic
program (QCQP): \begin{subequations} \label{eq_nn_5} 
\begin{align}
 & \min_{x\in\mathbb{R}^{n-1}}\ \ x^{*}{A}_{0}x+2b_{0}^{*}x+c_{0}\\
 & \ \ \mathrm{s.t.}\quad\hspace{3mm}x^{*}{A}_{k}x+2b_{k}^{*}x+c_{k}\leq0,\qquad k=1,\dots,m\label{QCQPcons}
\end{align}
\label{QCQP}\end{subequations} where $A_{k}\in\mathbb{S}^{n-1}$,
$b_{k}\in\mathbb{R}^{n-1}$ and $c_{k}\in\mathbb{R}$ for $k=0,\dots,m$.
This problem can be reformulated as 
\begin{equation}
\begin{aligned} & \min_{X\in\mathbb{S}^{n}}\quad\ \ \text{trace}\{M_{0}X\}\\
 & \ \ \ \mathrm{s.t.}\quad\;\;\,\text{trace}\{M_{k}X\}\leq0,\qquad k=1,\dots,m\\
 & \qquad\qquad X_{11}=1,\\
 & \qquad\qquad X\succeq0,\\
 & \qquad\qquad\text{rank}\{X\}=1
\end{aligned}
\label{eq_nn_6}
\end{equation}
where $X$ plays the role of 
\begin{equation}
\left[\begin{array}{cc}
1\\
x
\end{array}\right][1\quad x^{*}].\label{xfactor}
\end{equation}
and 
\begin{equation}
M_{k}=\left[\begin{array}{cc}
c_{k} & b_{k}^{*}\\
b_{k} & {A}_{k}
\end{array}\right],\qquad k=0,\dots,m
\end{equation}
Dropping the rank constraint from \eqref{eq_nn_6} leads to an SDP
relaxation of the QCQP problem~\eqref{eq_nn_4}. If the matrices
$M_{k}$'s are sparse, then the SDP relaxation is guaranteed to have
a low-rank solution. To enforce the low-rank solution to be rank-1,
one can use the idea described in~\eqref{eq_nn_2} and find a penalized
SDP approximation of \eqref{eq_nn_4} by first dropping the rank constraint
from \eqref{eq_nn_6} and then adding a penalty term similar to $\sum_{(i,j)\in\bar{\mathcal{E}}\backslash\mathcal{E}}t_{ij}\ {X}_{ij}$
to its objective function.

In an effort to study low-rank solutions of the SDP relaxation of
\eqref{eq_nn_4}, let $\mathcal{G}=(\mathcal{V},\mathcal{E})$ be
a graph with $n$ vertices such that $(i,j)\in\mathcal{G}$ if the
$(i,j)$ entry of at least one of the matrices $M_{0},M_{1},...,M_{m}$
is nonzero. The graph $\mathcal{G}$ captures the sparsity of the
optimization problem~\eqref{eq_nn_4}. Notice that those off-diagonal
entries of $X$ that correspond to non-existent edges of $\mathcal{G}$
play no direct role in the SDP relaxation. Let $\widehat{X}$ denote
an arbitrary solution of the SDP relaxation of \eqref{eq_nn_4}. It
can be observed that every solution to the low-rank positive-semidefinite
matrix completion problem \eqref{eq_nn_1} is a solution of the SDP
relaxation as well. Now, one can use the SDP problem \eqref{eq_nn_2}
to find low-rank feasible solutions of the SDP relaxation of the QCQP
problem. By carefully picking the graph $\bar{\mathcal{G}}$, one
can obtain a feasible solution of the SDP relaxation whose rank is
less than or equal to the treewidth of $\mathcal{G}$ plus 1 (this
number is expected to be small for sparse graphs) \cite{laurent2012new,Madani_SIAM_2017}.

\subsection{State Estimation}

Consider the problem of recovering the state of a nonlinear system
from noisy data. Without loss of generality, we only focus on the
quadratic case, where the goal is to find an unknown state/solution
$x\in\mathbb{R}^{n}$ for a system of quadratic equations of the form
\begin{align}
z_{r}={x}^{\ast}{M}_{r}x+\omega_{r},\qquad\forall r\in\{1,\ldots,m\}\label{feas}
\end{align}
where 
\begin{itemize}
\item $z_{1},\ldots,z_{m}\in\mathbb{R}$ are given measurements/specifications. 
\item Each of the parameters $\omega_{1},\ldots,\omega_{m}$ is an unknown
measurement noise with some known statistical information. 
\item ${M}_{1},\ldots,{M}_{m}$ are constant $n\times n$ matrices. 
\end{itemize}
Several algorithms in different contexts, such as signal processing,
have been proposed in the literature for solving special cases of
the above system of quadratic equations \cite{beck2013sparsity,beck2013nonlinear,chen2015solving,candes-2013,candes2015phase2,candes2013}.
These methods are often related to semidefinite programming. As an
example, consider the conic program \begin{subequations}
\begin{align}
 & \underset{\begin{subarray}{c}
\!{X}\in\mathbb{S}^{n}\\
\;\;{\nu}\in\mathbb{R}^{m}
\end{subarray}}{\min} &  & \text{trace}\{MX\}+\mu\left(\frac{|\nu_{1}|}{\sigma_{1}}\!+\!\cdots\!+\!\frac{|\nu_{m}|}{\sigma_{m}}\right)\label{sdpMnoisy_obj}\\
 & \ \ \text{s.t.} &  & \text{trace}\{M_{r}X\}+\nu_{r}=z_{r}, & \hspace{-15mm}r=1,\ldots,m\\
 &  &  & {X}\succeq0\label{sdpMpsdNoise}
\end{align}
\label{sdpMnoisy}\end{subequations} where $\mu>0$ is a sufficiently
large fixed parameter. The objective function of this problem has
two terms: one taking care of non-convexity and another one for estimating
the noise values. The matrix $M$ can be designed such that the solution
$X^{\text{opt}}$ of the above conic program and the unknown state
$x$ be related through the inequality: 
\begin{align}
\|{X}^{\mathrm{opt}}-\alpha{x}{x}^{*}\|\leq2\sqrt{\frac{\mu\times\|{\omega}\|_{1}\times\mathrm{trace}\{{X}^{\mathrm{opt}}\}}{\eta}}
\end{align}
where ${\omega}:=[\omega_{1}\ \cdots\ \omega_{m}]$, and $\alpha$
and $\eta$ are constants \cite{Madani_ES_CDC_2016}. This implies
that the distance between the solutions of the conic program and the
unknown state depends on the power of the noise. A slightly different
version of the above result holds even in the case where a modest
number of the values $\omega_{1},\ldots,\omega_{m}$ are arbitrarily
large corresponding to highly corrupted/manipulated data \cite{HICSS_Ramtin_2017}.
In that case, as long as the number of such measurements is not relatively
large, the solution of the conic program will not be affected. The
proposed technique has been successfully tested on the European Grid,
where more than 18,000 parameters were successfully estimated in \cite{Madani_ES_CDC_2016_2,Madani_ES_CDC_2016}.

\section{Convexification Techniques}

\label{sec:3}

In this section, we present convex relaxations for solving difficult
non-convex optimization problems. These non-convex problems include
integer programs, quadratically-constrained quadratic programs, and
more generally polynomial optimization. There are countless examples:
MAX CUT, MAX SAT, traveling salesman problem, pooling problem, angular
synchronization, optimal power flow, etc. Convex relaxations are crucial
when searching for integer solutions (e.g. via branch-and-bound);
when optimizing over real variables, they provide a way to find globally
optimal solutions, as opposed to local solutions. Our focus is on
hierarchies of relaxations that grow tighter and tighter towards the
original non-convex problem of interest. The common framework we consider
is that of optimizing a polynomial $f$ of $n$ variables constrained
by $m$ polynomial inequalities, i.e. 
\[
\inf_{x\in\mathbb{R}^{n}}~f(x)~~~\text{s.t.}~~~g_{i}(x)\geqslant0,~i=1,\hdots,m.
\]
We consider linear programming (LP) hierarchies, second order-conic
programming (SOCP) hierarchies, and semidefinite positive positive
(SDP) hierarchies. After briefly recalling some of the historical
contributions, we refer to recent developments that have taken place
over the last three years. They are aimed at making the hierarchies
more tractable. In particular, the recently proposed \textit{multi-ordered
Lasserre hierarchy} can solve a key industrial problem with thousands
of variables and constraints to global optimality. The Lasserre hierarchy
was previously limited to small-scale problems since it was introduced
17 years ago.

\subsection{LP hierarchies}

When solving integer programs, it is quite natural to consider convex
relaxations in the form of linear programs. The idea is to come up
with a polyhedral representation of the convex hull of the feasible
set so that all the vertices are integral. In that case, optimizing
over the representation can yield the desired integral solutions.
For an excellent reference on the various approaches, including those
of Gomori and Chvátal, see the recent book \cite{conforti2014}. It
is out of this desire to obtain nice representations that the Sherali-Adams
\cite{sherali1990} and Lovász-Schrijver \cite{lovasz1991} hierarchies
arose in 1990; they both provide tighter and tighter outer approximations
of the convex hull of the feasible set. In fact, after a finite number
of steps (which is known \textit{a priori}), their polyhedral representations
coincide with the convex hull of the integral feasible set. At each
iteration of the hierarchy, the representation of Sherali-Adams is
contained in the one of Lovász-Schrijver \cite{laurent-2003}. One
way to view these hierarchies is via lift-and-project. For instance,
the Sherali-Adams hierarchy can be viewed as taking products of the
constraints, which are redundant for the original problem, but strengthen
the linear relaxation. Interestingly, the well-foundedness of this
approach, and indeed the convergence of the Sherali-Adams hierarchy,
can be justified by the works of Krivine \cite{krivine1964,krivine-1964}
in 1964, as well as those of Cassier \cite{cassier1984} (in 1984)
and Handelman \cite{handelman1988} (in 1988). It was Lasserre \cite[Section 5.4]{lasserre-2010}
who recognized that, thanks to Krivine, one can generalize the approach
of Sherali-Adams to polynomial optimization. Global convergence is
ensured under some mild assumptions which can always be met when the
feasible set is compact. To do so, one needs to make some adjustements
to the modeling of the feasible set. These include normalizing the
constraints to be between zero and one. 
\begin{example}
\textit{Consider the following polynomial optimization problem taken
from} \cite[Example 5.5]{lasserre-2010}: 
\[
\inf_{x\in\mathbb{R}}~x(x-1)~~~\text{s.t.}~~~x\geqslant0~~\text{and}~~1-x\geqslant0
\]
\textit{Its optimal value is} $-1/4$\textit{. To obtain the second-order
LP relaxation, one can add the following redundant constraints:} 
\[
x^{2}\geqslant0,~~~x(1-x)\geqslant0,~~~(1-x)^{2}\geqslant0
\]
\textit{The lifted problem then reads:} 
\[
\inf_{y_{1},y_{2}\in\mathbb{R}}~y_{2}-y_{1}~~~\text{s.t.}~~~\left\{ \begin{array}{r}
y_{1}\geqslant0\\
1-y_{1}\geqslant0\\
y_{2}\geqslant0\\
y_{1}-y_{2}\geqslant0\\
1-2y_{1}+y_{2}\geqslant0
\end{array}\right.
\]
\textit{where} $y_{k}$ \textit{corresponds to} $x^{k}$ \textit{for
a positive integer} $k$\textit{. An optimal solution to this problem
is} $(y_{1},y_{2})=(1/2,0)$. \textit{We then obtain a lower bound
on the original problem equal to} $-1/2$\textit{. The third-order
LP relaxation is obtained by adding yet more redundant constraints:}
\[
x^{3}\geqslant0,~~~x^{2}(1-x)\geqslant0,~~~x(1-x)^{2}\geqslant0~~~(1-x)^{3}\geqslant0
\]
\textit{Now, the lifted problem reads:} 
\[
\inf_{y_{1},y_{2}\in\mathbb{R}}~y_{2}-y_{1}~~~\text{s.t.}~~~\left\{ \begin{array}{r}
y_{1}\geqslant0\\
1-y_{1}\geqslant0\\
y_{2}\geqslant0\\
y_{1}-y_{2}\geqslant0\\
1-2y_{1}+y_{2}\geqslant0\\
y_{3}\geqslant0\\
y_{2}-y_{3}\geqslant0\\
y_{1}-2y_{2}+y_{3}\geqslant0\\
1-3y_{1}+3y_{2}-y_{3}\geqslant0
\end{array}\right.
\]
\textit{An optimal solution is given by} $(y_{1},y_{2},y_{3})=(1/3,0,0)$\textit{,
yielding the lower bound of} $-1/3$\textit{. And so on and so forth.} 
\end{example}
While convergence of the Sherali-Adams hierarchy is preserved when
generalizing their approach to polynomial optimization, \textit{finite}
convergence is not preserved. In the words of Lasserre \cite[Section 5.4.2]{lasserre-2010}:
``\textit{Unfortunately, we next show that in general the LP-relaxations
cannot be exact, that is, the convergence is only asymptotic, not
finite.}'' This means that, while the global value can be approached
to abritrary accuracy, it may never be reached. As a consequence,
one cannot hope to extract global minimizers, i.e. those points that
satisfy all the constraints and whose evaluations are equal to the
global value.

We now turn our attention to some recent work on designing LP hierarchies
for polynomial optimization on the positive orthant, i.e. with the
constraints $x_{1}\geqslant0,\hdots x_{n}\geqslant0$. Invoking a
result of Póyla in 1928 \cite{polya1928}, it was recently proposed
in \cite{pena2017} to multiply the objective function by $(1+g_{1}(x)+\hdots+g_{m}(x))^{r}$
for some positive integer $r$, in addition to taking products of
the constraints (as in the aforementioned LP hierarchy). This work
comes as a result of generalizing the notion of copositivity from
quadratic forms to polynomial functions \cite{pena2015}.

We next discuss some numerical aspects of LP hierarchies for polynomial
optimization. It has been shown that multiplying constraints by one
another when applying lift-and-project to integer programs is very
efficient in practice. However, it can lead to numerical issues for
polynomial optimization. The reason for this is that the coefficients
in the polynomial constraints are not necessarily all of the same
order of magnitude. For example, in the univariate case, multiplying
$0.1x-2\geqslant0$ by itself yields $0.01x^{2}-0.2x+4\geqslant0$,
leading to coefficients ranging two orders of magnitude. This can
be challenging, even for state-of-the-art software in linear programming.
To date, there exists no way of scaling the coefficients of a polynomial
optimization problem so as to make them more or less homogenous. In
contrast, in integer programs, the constraints $x_{i}^{2}-x_{i}=0$
naturally have all coefficients equal to zero or one. As can be read
in \cite[page 7]{ahmadi-2017}: ``\textit{We remark that while there
have been other approaches to produce LP hierarchies for polynomial
optimization problems (e.g., based on the Krivine-Stengle certificates
of positivity} \cite{krivine-1964,lasserre-2010,stengle-1974}\textit{),
these LPs, though theoretically significant, are typically quite weak
in practice and often numerically ill-conditioned} \cite{toh-2017}.''
This leads us to discuss the LP hierarchy proposed in \cite{ahmadi-2014}
in 2014.

Viewed through the lenses of lift-and-project, the approach in \cite{ahmadi-2014}
avoids products of constraints and instead adds the redundant constraints
$x^{2\alpha}g_{i}(x)\geqslant0$ and $(x^{\alpha}\pm x^{\beta})^{2}g_{i}(x)\geqslant0$
where $x^{\alpha}=x_{1}^{\alpha_{1}}\cdots x_{n}^{\alpha_{n}}$ and
$\alpha,\beta\in\mathbb{N}^{n}$. By doing this for monomials of higher
and higher degree $\alpha_{1}+\hdots+\alpha_{n}$, one obtains an
LP hierarchy. The nice property of this hierarchy is that it avoids
the conditioning issues associated with the previously discussed hierarchies.
The authors also propose to multiply the objective by $(x_{1}^{2}+\hdots+x_{n}^{2})^{r}$
for some integer $r$ as a means to strenghen the hierarchy. Global
convergence can then be guaranteed (upon reformulation) when optimizing
a homogenous polynomial whose individual variables are raised only
to even degrees and are constrained to lie in the unit sphere. This
follows from \cite[Theorem 13]{ahmadi-2017}, a consequence of Pólya's
previously mentioned result \cite{polya1928}. 
We conclude by noting that the distinct LP hierarchies presented above
can be combined. For the interested reader, numerical experiments
can be found in \cite[Table 4.4]{kuang-2017}.

\subsection{SOCP hierarchies}

A natural way to provide hierarchies that are stronger than LP hierarchies
is to resort to conic optimization whose feasible set is not a polyhedra.
In fact, in their original paper, Lovász and Schrijver proposed strengthening
their LP hierarchy by adding a positive semidefinite constraint. We
will deal with SDP in the next section, and we now focus on SOCP for
which very efficient solvers exist. It was this practical consideration
which led the authors of \cite{ahmadi-2014} to restrain the cone
of sum-of-squares arising in the Lasserre hierarchy. By restricting
the number of terms inside the squares to two at most, i.e. by avoiding
squares such as the one crossed out in $\sigma(x_{1},x_{2})=x_{1}^{2}+(2x_{1}-x_{2}^{3})^{2}+\xcancel{(1-x_{1}+x_{2})^{2}}$,
one obtains SOCP constraints instead of SDP constraints. Numerical
experiments can found in \cite{kuang-2017bis}. A dual perspective
to this approach is to relax the semidefinite constraints in the Lasserre
hierarchy to all necessary SOCP constraints. This idea was independently
proposed in \cite{dan2015}, except that the authors of that work
do not relax the moment matrix. This ensures that the relaxation obtained
is at least as tight as the first-order Lasserre hiearchy. When applied
to find global minimizers to the optimal power flow problem \cite{carpentier-1962},
this provides reduced runtime in some instances. On other instances,
the hierarchy does not seem to globally converge, at least with the
limited computional power used in the experiments. This was elucidated
in \cite{josz2017} where it was shown that restricting sum-of-squares
to two terms at most does not preserve global convergence, even if
the polynomial optimization problem is convex. Interestingly, the
restriction on the sum-of-squares can be used to strenghten the LP
hierarchies described in the previous section \cite{kuang-2017}.
However, this does not affect the ill-conditioning associated with
the LP hierarchies, nor their asymptotic convergence (as opposed to
finite convergence).

\subsection{SDP hierarchies}

SDP hierarchies revolve around the notion of sum-of-squares, which
were first introduced in the context of optimization by Shor \cite{shor1987}
in 1987. Shor showed that globally minimizing a univariate polynomial
on the real line breaks down to a convex problem. It relies on the
fact that a univariate polynomial is nonnegative if and only if it
is a sum-of-squares of other polynomials. This is also true for bivariate
polynomials of degree four as well as multivariate quadratic polynomials.
But it is generally not true for other polynomials, as was shown by
Hilbert in 1888 \cite{hilbert1888}. This led Shor \cite{shor1998}
(see also \cite{ferrier1998}) to tackle the minimization of multivariate
polynomials by reformulating them using quadratic polynomials. Later,
Nesterov \cite{nesterov2000} provided a self-concordant barrier that
allows one to use efficient interior point algorithms to minimize
a univariate polynomial via sum-of-squares.

We turn our attention to the use of sum-of-squares in a more general
context, i.e. constrained optimization. Working on Markov chains where
one seeks invariant measures, Lasserre \cite{lasserre-2001} realized
that minimizing a polynomial function under polynomial constraints
can also be viewed as a problem where one seeks a measure. He showed
that a dual perspective to this approach consists in optimizing over
sum-of-squares. In order to justify the global convergence of his
approach, Lasserre used Putinar's Positivstellensatz \cite{putinar-1993}.
This result was discovered in 1993 and provided a crucial refinement
of Schmüdgen's Positivstellensatz \cite{schmudgen-1991} proven a
few years earlier. It was crucial because it enabled numerical computations,
leading to what is known today as the Lasserre hierarchy. 
\begin{example}
\textit{Consider the following polynomial optimization problem taken
from} \cite{josz2017}: 
\[
\inf_{x_{1},x_{2}\in\mathbb{R}}~x_{1}^{2}+x_{2}^{2}+2x_{1}x_{2}-4x_{1}-4x_{2}~~\text{s.t.}~~x_{1}^{2}+x_{2}^{2}=1
\]
\textit{Its optimal value is} $2-4\sqrt{2}$\textit{, which can be
found using sums-of-squares since:} 
\[
x_{1}^{2}+x_{2}^{2}+2x_{1}x_{2}-4x_{1}-4x_{2}~~~~-~~~~(2-4\sqrt{2})
\]
\[
=
\]
\[
(\sqrt{2}-1)(x_{1}-x_{2})^{2}+\sqrt{2}(-\sqrt{2}+x_{1}+x_{2})^{2}
\]
\[
+
\]
\[
2(\sqrt{2}-1)(1-x_{1}^{2}-x_{2}^{2})
\]
\textit{It can be seen from the above equation that when $(x_{1},x_{2})$
is feasible, the first line must be nonnegative, proving that $2-4\sqrt{2}$
is a lower bound. This corresponds to the first-order Lasserre hierarchy
since the polynomials inside the squares are of degree one at most.
To make the link with previous section, note that one can ask to restrict
the number of terms to two at most inside the squares. This allows
one to use second-order conic programming instead of semidefinite
programming, but does not preserve global convergence. The best bound
that can be obtained (i.e. $-4\sqrt{2}$) is given by the following
decomposition:} 
\[
x_{1}^{2}+x_{2}^{2}+2x_{1}x_{2}-4x_{1}-4x_{2}~~~~-~~~~(-4\sqrt{2})
\]
\[
=
\]
\[
\frac{\sqrt{2}}{2}(-\sqrt{2}+2x_{1})^{2}+\frac{\sqrt{2}}{2}(-\sqrt{2}+2x_{2})^{2}+(x_{1}+x_{2})^{2}
\]
\[
+
\]
\[
2\sqrt{2}(1-x_{1}^{2}-x_{2}^{2})
\]
\end{example}
Parallel to Lasserre's contribution, Parrilo \cite{parrilo-2000b}
pioneered the use of sum-of-squares for obtaining strong bounds on
the optimal solution of nonconvex problems (e.g. MAX CUT). He also
showed how they can be used for many important problems in systems
and control. These include Lyapunov analysis for control systems \cite{ahmadi2014}.
We do not dwell on these as they are outside the scope of this tutorial,
but we also mention a different view of optimal control via occupation
measures \cite{lasserre2008}. In contrast to Lasserre, Parrilo's
work \cite{parrilo-2003} panders to Stengle's Positivstellensatz
\cite{stengle-1974}, which is used for proving infeasibility of systems
of polynomial equations. This result can be seen as a generalization
of Farkas' Lemma which certifies the emptiness of a polyhedral set.

As discussed above, the Lasserre hierarchy provides a sequence of
semidefinite programs whose optimal values converge (monotonically)
towards the global value of a polynomial optimization problem. This
is true provided that the feasible set is compact, and that a bound
$R$ on the radius of the set is known, so that one can include a
redundant ball constraint $x_{1}^{2}+\cdots+x_{n}^{2}\leqslant R^{2}$.
When modeling the feasible set in this matter, there is also zero
duality gap in each semidefinite program \cite{josz-2015}. This is
a crucial property when using path following primal-dual interior
point methods, which are some of the most efficient approaches for
solving semidefinite programs.

In contrast to LP hierarchies which have only asymptotic convergence
in general, the Lasserre hierarchy has finite convergence generically.
This means that for a given abritary polynomial optimization problem,
finite convergence will almost surely hold. It is Nie \cite{nie-2014}
who proved this result, which had been observed in practice ever since
the Lasserre hierarchy had been introduced. He relied on theorems
of Marshall \cite{marshall2006,marshall2009} which attempted to answer
the question: when can a nonnegative polynomial have a sum-of-squares
decomposition? In the Positivstellensätze discussed above, the assumption
of positivity is made, which only guarantees asymptotic convergence.
The result of Nie marks a crucial difference with the LP hierarchies
because it means that in practice, one can solve non-convex problems
exactly via a convex relaxation, whereas with LP hierarchies one may
only approximate them. In fact, when finite convergence is reached,
the Lasserre hierarchy not only provides the global value, but also
finds global minimizers, i.e. points that satisfy all the constraints
and whose evaluations are the global value. This last feature illustrates
a nice synergy between advances in optimization and advances on the
theory of moments, which we next discuss.

When the Lasserre hierarchy was introduced, the theory of moments
lacked a result to guarantee when global solutions could be extracted.
At the time of Lasserre's original paper, the theory only applied
to bivariate polynomial optimization \cite[Theorem 1.6]{curto1991}.
With the success of the Lasserre hierarchy, there was a growing need
for more theory to be developed. This theory was developed a few years
later by Curto and Fialkow \cite[Theorem 1.1]{curto-2005}. They showed
that it is sufficient to check a rank condition in the Lasserre hierarchy
in order to extract global minimizers (and in fact, the number of
minimizers is equal to the rank).

Interestingly, the same situation occured when the complex Lasserre
hierarchy was recently introduced in \cite{joszmolzahn2017}. The
Lasserre hierarchy was generalized to complex numbers in order to
enhance its tractability when dealing with polynomial optimization
in complex numbers, i.e. 
\[
\inf_{z\in\mathbb{C}^{n}}~f(z,\bar{z})~~~\text{s.t.}~~~g_{i}(z,\bar{z})\geqslant0,~i=1,\hdots,m.
\]
where $f,g_{1},\hdots,g_{m}$ are real-valued complex polynomials
(e.g. $2|z_{1}|^{2}+(1+\textbf{i})z_{1}\bar{z}_{2}+(1-\textbf{i})\bar{z}_{1}z_{2}$,
where $\textbf{i}$ is the imaginary number). This framework is natural
for optimization problems with oscillatory phenoma, which are omnipresent
in physical systems (e.g. electric power systems, imaging science,
signal processing, automatic control, quantum mechanics). One way
of viewing the complex Lasserre hierarchy is that it restricts the
sums-of-squares in the original Lasserre hierarchy to \textit{Hermitian}
sums-of-squares. These are exponentially cheaper to compute yet preserve
global convergence, thanks to D'Angelo's and Putinar's Positivstellensatz
\cite{angelo-2008}. On the optimal power flow problem in electrical
engineering, they permit a speed-up factor of up to one order of magnitude
\cite[Table 1]{joszmolzahn2017}. 
\begin{example}
\textit{Consider the following complex polynomial optimization problem}
\[
\inf_{z\in\mathbb{C}}~~z+\bar{z}~~~\text{s.t.}~~~|z|^{2}=1
\]
\textit{whose optimal value is} $-2$\textit{. One way to solve this
problem would be to convert it into real numbers} $z=:x1+\textbf{i}x_{2}$\textit{,
i.e.} 
\[
\inf_{x_{1},x_{2}\in\mathbb{R}}~~2x_{1}~~~\text{s.t.}~~~x_{1}^{2}+x_{2}^{2}=1
\]
\textit{and to use sums-of-squares:} 
\[
2x_{1}~-~(-2)~~=~~1^{2}+(x_{1}+x_{2})^{2}~~+~~1\times(1-x_{1}^{2}-x_{2}^{2})
\]
\textit{But one could instead use Hermitian sums-of-squares which
are cheaper to compute:} 
\[
z+\bar{z}~~-~(-2)~~=~~~~~~~|1+z|^{2}~~~~~~+~~~~1\times(1-|z|^{2})
\]
\end{example}
When the complex Lasserre hierarchy was introduced, the theory of
moments lacked a result to guarantee when global solutions could be
extracted. This led its authors to generalize the work of Curto and
Fialkow using the notion of hyponormality in operator theory \cite[Theorem 5.1]{joszmolzahn2017}.
They found that in addition to rank conditions, some positive semidefinite
conditions must be met. Contrary to the rank conditions, these are
convex and can thus be added to the complex Lasserre hierarchy. In
\cite[Example 4.1]{joszmolzahn2017}, doing so reduces the rank from
3 to 1 and closes the relaxation gap.

The advent of the Lasserre hierarchy not only sparked progress in
the theory of moments, but also led to some notable results. In 2004,
the author of \cite{schweighofer2004} derived an upper bound on the
order of the Lasserre hierarchy needed to obtain a desired relaxation
bound to a polynomial optimization problem. This upper bound depends
on three factors: 1) a certain description of the feasible set, 2)
the degree of the polynomial objective function, and 3) how close
the objective function is to reaching the relaxation bound on the
feasible set. It must be noted that this upper bound is difficult
to compute in practice. As of today, it is therefore not possible
to know ahead of time how far in the hierarchy one needs to go in
order to solve a given instance of polynomial optimization. Along
these lines, nothing is yet known about the speed of convergence of
the bounds generated by the Lasserre hierarchy. In practice, they
generally reach the global value in a few iterations. Somewhat paradoxically,
there are some nice results \cite{deklerk2017} on the speed of convergence
of SDP hierarchies of upper bounds \cite{lasserre2011}, although
their converge is slow in practice.

Another notable result is contained in \cite{blekherman2006}. As
discussed previously, the discrepancy between nonnegative polynomials
and sum-of-squares was noticed towards the end of the nineteenth century.
In some applications of sum-of-squares, the nonnegativity of a function
on $\mathbb{R}^{n}$ is replaced by requiring it to be a sum-of-squares.
The result in \cite{blekherman2006} quantifies how small the cone
of sums-of-squares is with respect to the cone of nonnegative polynomials.
It is perhaps the title of the paper that best sums up the finding:
``\textit{There are significantly more nonnegative polynomials than
sums of squares}''.

Having discussed several theoretical aspects of the Lasserre hierarchy,
we now turn our attention to practical considerations. In order to
make the Lasserre hierarchy tractable, it is crucial to exploit the
problem structure. We have already gotten a flavor of this with the
complex hierarchy, which exploits the complex structure of physical
problems with oscillatory phenoma (electricity, light, etc.). In the
following, some key results on sparsity and symmetry are highlighted.

In order to exploit sparsity in the Lasserre hierarchy, it was proposed
to use chordal sparsity in sums-of-squares in \cite{waki-2006}. We
briefly explain this approach. Each constraint in a polynomial optimization
problem is associated a sum-of-squares in the Lasserre hierarchy.
In fact, each sum-of-squares can be interpreted as a generalized Lagrange
multiplier \cite[Theorem 5.12]{lasserre-2010}. If a constraint only
depends on a few variables, say $x_{1},x_{3},x_{20}$ among $x_{1},\hdots,x_{100}$,
it seems naturally that the associated sum-of-squares should depend
only on the variables $x_{1},x_{3},x_{20}$, or some slightly larger
set of variables. This was made possible by the work in \cite{waki-2006}.
To do so, the authors consider the correlative sparsity pattern of
the polynomial optimization problem. It can be viewed as a graph where
the nodes are the variables and the edges signify a coupling of the
variables. Taking a chordal extension of this graph and computing
the maximal cliques, one can restraint the variables appearing in
the sum-of-squares to belong to these cliques. This provides a more
tractable hierarchy of relaxations while preserving global convergence,
as was shown by Lasserre \cite[Theorems 2.28 and 4.7]{lasserre-2010}.
What Lasserre proved was a sparse version of Punitar's Positivstellensatz.
This sparse version has two applications that we next discuss.

One application is to the bounded sum-of-squares hierarchy (BSOS)
\cite{toh-2017}. The idea of this SDP hierarchy is to fix the size
of the SDP constraint as the order of the hierarchy increases. The
size of the SDP constraint can be set by the user. The hierarchy builds
on the LP hierarchy based on Krivine's Positivstellensatz discussed
in the section on LP hierarchies. As the order of the hierarchy increases,
the number of LP constraints augments. These are the LP constraints
that arise when multiplying constraints by one another. A sparse BSOS
\cite{weisser-2017} is possible, thanks to the sparse version of
Putinar's Positivstellensatz. Global convergence is guaranteed by
\cite[Theorem 1]{weisser-2017}, but the ill-conditioning associated
with LP hierarchies is inherited.

Another applicaton of the sparse version of Punitar's Positivstellensatz
is the \textit{multi-ordered} Lasserre hierarchy \cite{mh_sparse_msdp,joszmolzahn2017}.
It is based on two ideas: 1) to use a different relaxation order for
each constraint, and 2) to iteratively seek a closest measure to the
truncated moment data until a measure matches the truncated data.
Global convergence is a consequence of the aforementioned sparse Positivstellensatz.
\begin{example}
\textit{The multi-ordered Lasserre hierarchy can solve a key industrial
problem of the twentieth century to global optimality on instances
of polynomial optimization with up to} 4,500 \textit{variables and}
14,500 \textit{constraints} (\textit{see table below and} \cite{joszmolzahn2017})\textit{.
The relaxation order is typically augmented at a hundred or so constraints
before reaching global optimality. The test cases correspond to the
highly non-convex optimal flow problem, and in particular to instances
of the European high-voltage synchronous electricity network comprising
data from} 23 \textit{different countries} (\textit{available at}
\cite{josz-2016,nesta}). 
\begin{table}[!ht]
\centering %
\begin{tabular}{|l|c|c|c|c|}
\hline 
\multicolumn{1}{|c|}{\textbf{Case}} & \multicolumn{1}{c|}{\textbf{Number of}} & \multicolumn{1}{c|}{\textbf{Const-}} & \multicolumn{2}{c|}{\textbf{Multi-ordered Lasserre}}\tabularnewline
\cline{4-5} 
\multicolumn{1}{|c|}{\textbf{Name}} & \multicolumn{1}{c|}{\textbf{Variables}} & \multicolumn{1}{c|}{\textbf{raints}} & \hspace*{1em}\textbf{Global val.}\hspace*{1em}  & \textbf{Time (s.)} \tabularnewline
\hline 
case57Q  & \hphantom{1,}114  & \hphantom{11,}192  & \hphantom{11}7,352  & \hphantom{1,11}3.4 \tabularnewline
case57L  & \hphantom{1,}114  & \hphantom{11,}352  & \hphantom{1}43,984  & \hphantom{1,11}1.4 \tabularnewline
case118Q  & \hphantom{1,}236  & \hphantom{11,}516  & \hphantom{1}81,515  & \hphantom{1,1}15.7 \tabularnewline
case118L  & \hphantom{1,}236  & \hphantom{11,}888  & 134,907  & \hphantom{1,1}10.5 \tabularnewline
case300  & \hphantom{1,}600  & \hphantom{1}1,107  & 720,040  & \hphantom{1,11}7.2 \tabularnewline
nesta\_case24  & \hphantom{1,1}48  & \hphantom{11,}526  & \hphantom{11}6,421  & \hphantom{1,}246.1\tabularnewline
nesta\_case30  & \hphantom{1,1}60  & \hphantom{11,}272  & \hphantom{111,}372  & \hphantom{1,}302.7 \tabularnewline
nesta\_case73  & \hphantom{1,}146  & \hphantom{1}1,605  & \hphantom{1}20,125  & \hphantom{1,}506.9 \tabularnewline
PL-2383wp  & 4,354  & 12,844  & \hphantom{1}24,990  & \hphantom{1,}583.4 \tabularnewline
PL-2746wop  & 4,378  & 13,953  & \hphantom{1}19,210  & 2,662.4 \tabularnewline
PL-3012wp  & 4,584  & 14,455  & \hphantom{1}27,642  & \hphantom{1,}318.7 \tabularnewline
PL-3120sp  & 4,628  & 13,948  & \hphantom{1}21,512  & \hphantom{1,}386.6 \tabularnewline
PEGASE-1354  & 1,966  & \hphantom{1}6,444  & \hphantom{1}74,043  & \hphantom{1,}406.9 \tabularnewline
PEGASE-2869  & 4,240  & 12,804  & 133,944  & \hphantom{1,}921.3 \tabularnewline
\cline{1-4} 
\end{tabular}\label{tab:results} 
\end{table}

\textit{MOSEK's interior point software for semidefinite programming
is used in numerical experiments. Better runtimes (by up to an order
of magnitude) and higher precision can be obtained with the complex
Lasserre hierarchy mentioned above.} 
\end{example}
We finish by discussing symmetry. In the presence of symmetries in
a polynomial optimization problem, the authors of \cite{riener-2013}
proposed to seek an invariant measure when deploying the Lasserre
hierarchy. This reduces the computational burden in the semidefinite
relaxations. It was shown recently that in the presence of commonly
encountered symmetries, one actually gets a \textit{block diagonal}
Lasserre hierarchy \cite[Section 7]{joszmolzahn2017}. 
\begin{example}
\textit{Consider the following polynomial optimization problem} 
\[
\inf_{x_{1},x_{2}\in\mathbb{R}}~~2x_{1}x_{2}~~~\text{s.t.}~~~x_{1}^{2}+x_{2}^{2}=1
\]
\textit{whose optimal value is} $-1$\textit{. If} $(x_{1},x_{2})$
\textit{is feasible, then so is} $(-x_{1},-x_{2})$\textit{. As a
result, one may seek even sums-of-squares:} 
\[
2x_{1}x_{2}~-~(-1)~~=~~(x_{1}+x_{2})^{2}~+~1\times(1-x_{1}^{2}-x_{2}^{2})
\]
\textit{The sum-of-squares is even because it only has monomials of
even degree:} 
\[
(x_{1}+x_{2})^{2}~=~x_{1}^{2}+2x_{1}x_{2}+x_{2}^{2}
\]
\textit{Searching for even sums-of-squares reduces to block diagonal
SDP's} (\textit{see} \cite[Section 7]{joszmolzahn2017} \textit{for
explanations}). 
\end{example}
The above approaches for making hierarchies more tractable lead to
convex models that are more amenable for off-the-shelf solvers. A
recent approach \cite{zheng2017} proposes to adapt the algorithm
in the solver to exploit the structure brought about by sums-of-squares.
The next section deals with numerical algorithms for general conic
optimization problems.

\section{Numerical Algorithms}

\label{sec:4}

\global\long\def\R{\mathbb{R}}
 \global\long\def\S{\mathbb{S}}
 \global\long\def\cone{\mathcal{K}}
 \global\long\def\feas{\mathcal{F}}
 \global\long\def\R{\mathbb{R}}
 \global\long\def\C{\mathbb{C}}
 \global\long\def\S{\mathbb{S}}
 \global\long\def\A{\mathbf{A}}
 \global\long\def\B{\mathbf{B}}
 \global\long\def\H{\mathbf{H}}
 \global\long\def\P{\mathbf{P}}
 \global\long\def\D{\mathbf{D}}
 \global\long\def\tr{\mathrm{trace}}
 \global\long\def\vector{\mathrm{vec}\,}
 \global\long\def\plane{\mathcal{H}}
 \global\long\def\proj{\mathrm{proj}}
 \global\long\def\prox{\mathrm{prox}}
 \global\long\def\diag{\mathrm{diag}}
 \global\long\def\U{\mathbf{U}}
 \global\long\def\rank{\mathrm{rank}}
 \global\long\def\Sch{\mathbf{S}}
 \global\long\def\Lag{\mathscr{L}}
 \global\long\def\opt{\mathrm{opt}}

The previous section described \emph{convexification} techniques that
relax hard, nonconvex problems into a handful of standard class of
convex optimization problems, all of which can be approximated to
arbitrary accuracy in polynomial time using the ellipsoid algorithm~\cite{nemirovskii1983problem,shor1985minimization}.
Whenever the relaxation is tight, a solution to the original nonconvex
problem can be recovered after solving the convexified problem. Accordingly,
convexification establishes the original problem to be tractable or
``easy to solve'', at least in a theoretical sense. This approach
was used as early as 1980 by Gröteschel, Lovász and Schijver~\cite{grotschel1981ellipsoid}
to develop polynomial-time algorithms for combinatorial optimization.

However, the practical usefulness of convexification was less clear
at the time of its development. The ellipsoid method was notoriously
slow in practice, so specialized algorithms had to be used to solve
the resulting convexified problems. The fastest was the simplex method
for the solution of linear programs (LPs), but LP convexifications
are rarely tight. Conversely, semidefinite program (SDP) convexifications
are often exact, but SDPs were particularly difficult to solve, even
in very small instances (see~\cite{boyd1994linear,vandenberghe1996semidefinite}
for the historial context). As a whole, convexification remained mostly
of theoretical interest.

In the 1990s, advancements in numerical algorithms overhauled this
landscape. The interior-point method—originally developed as a practical
but rigorous algorithm for LPs by Karmarkar~\cite{karmarkar1984new}—was
extended to SDPs by Nesterov and Nemirovsky~\cite[Chapter 4]{nesterov1994interior}
and Alizadeh~\cite{alizadeh1995interior}. In fact, this line of
work showed interior-point methods to be \emph{particularly suitable}
for SDPs, generalizing and unifying the much of the previous framework
developed for LPs. For control theorists, the ability to solve SDPs
from convexification had a profound impact, giving rise to the disciplines
of LMI control~\cite{boyd1994linear} and polynomial control~\cite{parrilo2000structured}. 

Today, the growth in the size of SDPs has outpaced the ability of
general-purpose interior-point methods to solve them, fueled in a
large part by the application of convexification techniques to control
and machine learning applications. First-order methods have become
popular, because they have very low per-iteration costs that can often
be custom-tailored to exploit problem structure in a specific application.
On the other hand, these methods typically require considerably more
iterations to converge to reasonable accuracy. Ultimately, the most
effective algorithms for the solution of large-scale SDPs are those
that combine the convergence guarantees of interior-point methods
with the ability of first-order methods to exploit problem structure. 

This section reviews in detail three numerical algorithms for SDPs.
First, we describe the theory of interior-point methods in Section~\ref{subsec:Interior-point-method},
presenting them as a general-purpose algorithm for solving SDPs in
polynomial time. Next, we describe a popular first-order method in
Section~\ref{subsec:ADMM} known as ADMM, and explain how it is able
to reduce computational cost by exploiting sparsity. In Section~\ref{subsec:mIPM},
we describe a modified interior-point method for low-rank SDPs that
exploits sparsity like ADMM while also enjoying the strong convergence
guarantees of interior-point methods. Finally, we briefly review other
structure-exploiting algorithms in Section~\ref{subsec:other}.

\subsection{Problem description}

In order to simplify our presentation, we will focus our efforts on
the standard form semidefinite program 
\begin{align}
X^{\opt}= & \text{ minimize } & C\bullet X & \tag{SDP}\label{eq:SDP}\\
 & \text{ subject to } & A_{i}\bullet X & =b_{i}\quad\forall i\in\{1,\ldots,m\}\nonumber \\
 &  & X & \succeq0,\nonumber 
\end{align}
over the data $C,A_{1},\ldots,A_{m}\in\S^{n}$, $b\in\R^{m}$, and
its Lagrangian dual 
\begin{align}
\{y^{\opt},S^{\opt}\}= & \text{ maximize } & b^{T}y & \tag{SDD}\label{eq:SDD}\\
 & \text{ subject to } & \sum_{i=1}^{m}y_{i}A_{i}+S & =C\nonumber \\
 &  & S & \succeq0.\nonumber 
\end{align}
Here, $X\succeq0$ indicates that $X$ is positive semidefinite, and
$A_{i}\bullet X=\tr\{A_{i}X\}$ is the usual matrix inner product.
In case of nonunique solutions, we use $\{X^{\opt},y^{\opt},S^{\opt}\}$
to refer to the \emph{analytic center} of the solution set. We make
the following nondegeneracy assumptions.

\begin{assume}[Linear independence]\label{ass:lin}The matrix $\A=[\vector A_{1},\ldots,\vector A_{m}]$
has full column-rank, meaning the matrix $\A^{T}\A$ is invertible.\end{assume}

\begin{assume}[Slater's Condition]\label{ass:slater}There exist
$y$ and positive definite $X$ and $S$ such that $A_{i}\bullet X=b_{i}$
holds for all $i$, and $\sum_{i}y_{i}A_{i}+S=C$.\end{assume}

These are generic properties of SDPs, and are satisfied by almost
all instances~\cite{alizadeh1997complementarity}. Linear independence
implies that the number of constraints $m$ cannot exceed the number
of degrees of freedom $\frac{1}{2}n(n+1)$. Slater's condition is
commonly satisfied by embedding (\ref{eq:SDP}) and (\ref{eq:SDD})
within a slightly larger problem using the homogenous self-dual embedding
technique~\cite{ye1994homogeneous}.

All of our algorithms and associated complexity bounds can be generalized
in a straightforward manner to conic programs posed on the Cartesian
product of many semidefinite cones $\cone=\S_{+}^{n_{1}}\times\S_{+}^{n_{2}}\times\cdots\times\S_{+}^{n_{\ell}}$,
as in 
\begin{align}
 & \text{ minimize } & \sum_{j=1}^{\ell}C_{j}\bullet X_{j}\label{eq:SDP-1}\\
 & \text{ subject to } & \sum_{j=1}^{\ell}A_{i,j}\bullet X_{j} & =b_{i}\quad\forall i\in\{1,\ldots,m\}\nonumber \\
 &  & X_{j} & \succeq0\quad\forall j\in\{1,\ldots,\ell\}\nonumber 
\end{align}
and 
\begin{align}
 & \text{ maximize } & b^{T}y\label{eq:SDD-1}\\
 & \text{ subject to } & \sum_{i=1}^{m}y_{i}A_{i,j}+S_{j} & =C_{j}\quad\forall j\in\{1,\ldots,\ell\}\nonumber \\
 &  & S_{j} & \succeq0\quad\forall j\in\{1,\ldots,\ell\}.\nonumber 
\end{align}
We will leave the specific details as an exercise for the reader.
Note that this generalization includes linear programs (LPs), since
the positive orthant is just the Cartesian product of many size-1
semidefinite cones, as in $\R_{+}^{n}=\S_{+}^{1}\times\cdots\times\S_{+}^{1}$.

At least in principle, our algorithms also generalize to second-order
cone programs (SOCPs) by converting them into SDPs, as in
\[
\|u\|_{2}\le u_{0}\iff\begin{bmatrix}u_{0} & u^{T}\\
u & u_{0}I
\end{bmatrix}\succeq0.
\]
However, as demonstrated by Lobo~et~al.~\cite{lobo1998applications},
considerable efficiency can be gained by treating SOCPs as its own
distinct class of conic problems.

\subsection{\label{subsec:Interior-point-method}Interior-point methods}

The original interior-point methods were inspired by the logarithmic
barrier method, which replaces each inequality constraint of the form
$c(x)\ge0$ is by a logarithmic penalty term $-\mu\log c(x)$ that
is well-defined at \emph{interior-points} where $c(x)>0$, but becomes
unbounded from above as $x$ approaches the boundary where $c(x)=0$.
(This behavior constitutes an infinite \emph{barrier} that restricts
$x$ to lie within the feasible region where $c(x)>0$.)

Consider applying this strategy to (\ref{eq:SDP}) and (\ref{eq:SDD}).
If we intepret the semidefinite condition $X\succeq0$ as a set of
eigenvalue constraints 
\[
\lambda_{j}(X)\ge0\text{ for all }j\in\{1,\ldots,n\},
\]
then the resulting logarithmic barrier is none other than the log-determinant
penalty for determinant maximization (see~\cite{vandenberghe1998determinant}
and the references therein) 
\[
-\sum_{j=1}^{n}\mu\log\lambda_{j}(X)=-\mu\log\prod_{j=1}^{n}\lambda_{j}(X)=-\mu\log\det X.
\]
Substituting the penalty in place of the constraints $X\succeq0$
and $S\succeq0$ results in a sequence of unconstrained problems 
\begin{align}
X_{\mu}=\text{ minimize } & C\bullet X-\mu\log\det X\tag{SDP\ensuremath{\mu}}\label{eq:PotPrim}\\
\text{subject to } & A_{i}\bullet X=b_{i}\;\forall i\in\{1,\ldots,m\},\nonumber 
\end{align}
and 
\begin{align}
\{y_{\mu},S_{\mu}\}=\text{ maximize } & b^{T}y+\mu\log\det S\tag{SDD\ensuremath{\mu}}\label{eq:PotDual}\\
\text{subject to } & \sum_{i=1}^{m}y_{i}A_{i}+S=C,\nonumber 
\end{align}
which can be shown to be primal-dual pairs (up to a constant offset).

After converting inequality constraints into logarithmic penalties,
the barrier method repeatedly solves the resulting unconstrained problem
using progressively smaller values of $\mu$, each time reusing the
most recent solution as the starting point for the next minimization.
Applying this sequential strategy to solve (\ref{eq:PotPrim}) using
Newton's method yields a \emph{primal-scaled} interior-point method;
doing the same for (\ref{eq:PotDual}) results in a \emph{dual-scaled}
interior-point method. It is a seminal result of Nesterov and Nemirovski~\cite{nesterov1994interior}
that either interior-point methods converge to an approximate solution
accurate to $L$ digits after at most $O(nL)$ Newton iterations.
(This can be further reduced to $O(\sqrt{n}L)$ Newton iterations
by limiting the rate at which $\mu$ is reduced.) In practice, convergence
almost never occurs in more than tens of iterations.

In finite precision, the primal-scaled and dual-scaled interior-point
methods can suffer from severe accuracy and robustness issues; these
are the same reasons that had originally caused the barrier method
to fall out favor in the 1970s. Today, the most robust and accurate
interior-point methods are \emph{primal-dual}, and simultaneously
solve (\ref{eq:PotPrim}) and (\ref{eq:PotDual}) through their joint
Karush–Kuhn–Tucker (KKT) optimality conditions\begin{subequations}\label{eq:central-path}
\begin{align*}
A_{i}\bullet X_{\mu} & =b_{i}\quad\forall i\in\{1,\ldots,m\},\\
\sum_{i=1}^{m}y_{i}^{\mu}A_{i}+S_{\mu} & =C,\\
X_{\mu}S_{\mu} & =\mu I.
\end{align*}
\end{subequations}Here, the barrier parameter $\mu>0$ controls the
duality gap between the point $X_{\mu}$ in (\ref{eq:SDP}) and the
point $\{y_{\mu},S_{\mu}\}$ in (\ref{eq:SDD}), as in 
\[
n\mu=X_{\mu}\bullet S_{\mu}=C\bullet X_{\mu}-b^{T}y_{\mu}.
\]
In other words, the candidate solutions $X_{\mu}$ and $\{y_{\mu},S_{\mu}\}$
are suboptimal for (\ref{eq:SDP}) and (\ref{eq:SDD}) respectively
by an absolute figure no worse than $n\mu$, 
\begin{align*}
C\bullet X^{\star} & \le C\bullet X_{\mu}\le C\bullet X^{\star}+n\mu,\\
b^{T}y^{\star}-n\mu & \le b^{T}y_{\mu}\le b^{T}y^{\star}.
\end{align*}
The solutions $\{X_{\mu},y_{\mu},S_{\mu}\}$ for different values
of $\mu$ trace a trajectory in the feasible region that approaches
$\{X^{\opt},y^{\opt},S^{\opt}\}$ as $\mu\to0^{+}$, known as the
\emph{central path}. Modern primal-dual interior-point methods for
SDPs like SeDuMi, SDPT3, and MOSEK use Newton's method to solve the
KKT equations (\ref{eq:central-path}), while keeping each iterate
within a \emph{wide neighborhood} of the central path 
\[
\mathcal{N}_{\infty}^{-}(\gamma)\triangleq\left\{ \{X,y,S\}\in\feas:\lambda_{\min}(XS)\ge\frac{\gamma}{n}X\bullet S\right\} ,
\]
where $\feas$ denotes the feasible region 
\[
\feas\triangleq\left\{ \{X,y,S\}:\begin{array}{rl}
A_{i}\bullet X & =b_{i}\;\forall i,\\
\sum_{i}y_{i}A_{i}+S & =C,\\
X,S & \text{ pos. def.}
\end{array}\right\} .
\]
Here, $\gamma\in(0,1)$ quantifies the ``size'' of the neighborhood,
and is typically chosen with an aggressive value like $10^{-3}$.
The resulting algorithm is guarantee to converge to an approximate
solution accurate to $L$ digits after at most $O(nL)$ Newton iterations.
(This can be further reduced to $O(\sqrt{n}L)$ Newton iterations
by adopting a narrow neighborhood.) In practice, convergence almost
always occurs with 30-50 iterations.

\subsubsection{Complexity}

All interior-point methods converge to $L$ accurate digits in between
$O(\sqrt{n}L)$ and $O(nL)$ Newton iterations, and practical implementations
almost always occurs with tens of iterations. Accordingly, the complexity
of solving (\ref{eq:SDP}) and (\ref{eq:SDD}) using an interior-point
method is—up to a small multiplicative factor—the same as the cost
of solving the associated Newton subproblem 
\begin{align}
\text{ maximize } & b^{T}y-\frac{1}{2}\|W^{\frac{1}{2}}(S-Z)W^{\frac{1}{2}}\|_{F}^{2}\label{eq:Newt}\\
\text{subject to } & \sum_{i=1}^{m}y_{i}A_{i}+S=C,\nonumber 
\end{align}
in which $W,Z\in\S_{++}^{n}$ are used by the algorithm to approximate
the log-det penalty\footnote{Here, we assume that primal-, dual-, or Nesterov–Todd primal-dual
scalings are used. The less common H..K..M and AHO primal-dual scalings
have a sightly different version of (\ref{eq:Newt}); see~\cite{todd1998nesterov}
for a comparison.}. The positive definite matrix $W$ is known as the \emph{scaling
matrix}, and is always fully-dense.

The standard approach found in the vast majority of interior-point
solvers is form and solve the \emph{Hessian equation} (also known
as the Schur complement equation), obtained by substituting $S=C-\sum_{i=1}^{m}y_{i}A_{i}$
into the objective (\ref{eq:Newt}) and taking first-order optimality
conditions: 
\begin{equation}
A_{i}\bullet\left[W\left(\sum_{j=1}^{m}y_{j}A_{j}\right)W\right]=\underbrace{b_{i}+A_{i}\bullet W(C-Z)W}_{r_{i}}\label{eq:hesseqn2}
\end{equation}
for all $i\in\{1,\ldots,m\}$. Once $y$ is computed, the variables
$S=C-\sum_{i=1}^{m}y_{i}A_{i}$ and $X=W(Z-S)W$ are easily recovered.
Vectorizing the matrix variables allows (\ref{eq:hesseqn2}) to be
compactly written as 
\begin{equation}
\underbrace{[\A^{T}(W\otimes W)\A]}_{\H}y=r\label{eq:hesseqn1}
\end{equation}
where $\A=[\vector A_{1},\ldots,\vector A_{m}]$. It is common to
solve (\ref{eq:hesseqn1}) by forming the Hessian matrix $\H$ explicitly
and factoring it using dense Cholesky factorization, in $O(n^{3}m+n^{2}m^{2}+m^{3})$
time and $\Theta(m^{2}+n^{2})$ memory. The overall interior-point
method then has complexity between $\sim n^{3}$ and $\sim n^{6}$
time, requiring $\sim m^{2}$ memory.

\subsubsection{Bibliography}

The modern study of interior-point methods was initiated by Karmarkar~\cite{Karmarkar1984}
and their extension to SDPs was due to Nesterov and Nemirovsky~\cite[Chapter 4]{nesterov1994interior}
and Alizadeh~\cite{alizadeh1995interior}. Earlier algorithms were
essentially the same as the barrier methods from the 1960s; M.~Wright~\cite{wright2005interior}
gives an overview of this historical context. The effectiveness of
these methods was explained by the seminal work of Nesterov and Nemirovski~\cite{nesterov1994interior}
on self-concordant barrier methods. For an accessible introduction
to these classical results, see Boyd and Vandenberghe~\cite[Chapters 9 and 11]{boyd2004convex}.
The development of primal-dual interior-point methods began with Kojima~et~al.~\cite{kojima1989algorithm,kojima1991unified},
and was eventually extended to semidefinite programming and second-order
cone programming in a unified way by Nesterov and Todd~\cite{nesterov1997self,nesterov1998primal}.
For a survey on primal-dual interior-point methods for SDPs, see Sturm~\cite{sturm2002implementation}.
Today, the best interior-point solvers for SDPs are SeDuMi, SDPT3,
and MOSEK; the interested reader is referred to~\cite{sturm2002implementation,tutuncu2003solving,andersen2000mosek}
for their implementation details.

\subsection{\label{subsec:ADMM}ADMM}

One of the most successful first-order methods for SDPs has been ADMM.
Part of its appeal is that it is simple and easy to implement at a
large scale, and that convergence is guaranteed under very mild assumptions.
Furthermore, the algorithm is often ``lucky'': for many large-scale
SDPs, it converges at a \emph{linear rate}—like an interior-point
method—to 6+ digits of accuracy in just a few hundred iterations~\cite{wen2010alternating}.
However, the worst-case behavior is regularly attained, particularly
for SDPs that arise from convexification; thousands of iterations
are required to obtain solutions of only modest accuracy~\cite{madani2015admm,zheng2016fast}.

ADMM, or the alternating direction method of multipliers, is closely
related to the augmented Lagrangian method, a popular optimization
algorithm in the 1970s for constrained optimization that was historically
known as ``the method of multipliers''~\cite{hestenes1969multiplier,powell1969method}.
Both methods begin by augmenting the dual problem (\ref{eq:SDD}),
written here as a minimization 
\begin{align}
 & \text{minimize } & -b^{T}y\label{eq:dual_prob}\\
 & \text{subject to } & \sum_{i=1}^{m}y_{i}A_{i}+S & =C\nonumber \\
 &  & S & \succeq0\nonumber 
\end{align}
with a quadratic penalty term that does not affect the minimum nor
the minimizer 
\begin{align}
 & \text{minimize } & -b^{T}y+\frac{t}{2}\left\Vert \sum_{i=1}^{m}y_{i}A_{i}+S-C\right\Vert _{F}^{2}\label{eq:aug_dual_prob}\\
 & \text{subject to } & \sum_{i=1}^{m}y_{i}A_{i}+S=C\nonumber \\
 &  & S\succeq0.\nonumber 
\end{align}
However, the quadratic term makes the problem (\ref{eq:aug_dual_prob})
strongly convex (for $t>0$ and under Assumption~\ref{ass:lin}).
The convex conjugate of a strongly convex function is Lipschitz smooth
(see e.g.~\cite{nesterov2005smooth}), and this mean that the Lagrangian
dual, written 
\begin{equation}
\underset{X\succeq0}{\text{maximize }}\left\{ \min_{y,S\succeq0}\Lag_{t}(X,y,S)\right\} ,\label{eq:aug_prim_prob}
\end{equation}
where $\Lag_{t}$ is the augmented Lagrangian 
\begin{multline*}
\Lag_{t}(X,y,S)=-b^{T}y+X\bullet\left(\sum_{i=1}^{m}y_{i}A_{i}+S-C\right)\\
+\frac{t}{2}\left\Vert \sum_{i=1}^{m}y_{i}A_{i}+S-C\right\Vert _{F}^{2},
\end{multline*}
is differentiable with a Lipschitz-continuous gradient. Essentially,
adding a quadratic regularization to the dual problem (\ref{eq:dual_prob})
smoothes the corresponding primal problem (\ref{eq:aug_prim_prob}),
thereby allowing a gradient-based optimization algorithm to be effectively
used for its solution.

The augmented Lagrangian algorithm is derived by applying projected
gradient ascent to the maximization (\ref{eq:aug_prim_prob}) while
setting the step-size to exactly $t$, as in 
\[
X^{k+1}=\left[X^{k}+t\nabla_{X}\left\{ \min_{y,S\succeq0}\Lag_{t}(X^{k},y,S)\right\} \right]_{+},
\]
where $[W]_{+}$ denotes the projection of the matrix $W$ onto the
semidefinite cone $\S_{+}^{n}$, as in 
\[
[W]_{+}=\arg\min_{Z\succeq0}\|W-Z\|_{F}^{2}.
\]
Some algebra shows that the gradient term can be evaluated by solving
the inner minimization problem, and that the special step-size of
$t$ guarantees $X^{k}\succeq0$ so long as $X^{0}\succeq0$. Substituting
these two simplifications yields the classic form of the augmented
Lagrangian sequence 
\begin{align*}
\{y^{k+1},S^{k+1}\} & =\arg\min_{y,S\succeq0}\Lag_{t}(X^{k},y,S)\\
X^{k+1} & =X^{k}+t\left(\sum_{i=1}^{m}y_{i}^{k+1}A_{i}+S^{k+1}-C\right).
\end{align*}
The iterates converge to the solutions of (\ref{eq:SDP}) and (\ref{eq:SDD})
for all fixed $t>0$, and that the convergence rate is super-linear
if $t$ is allowed to increase after every iteration~\cite{rockafellar1976monotone}.
In practice, convergence to high accuracy is achieved in tens of iterations
by picking a very large value of $t$.

A key difficulty of the augmented Lagrangian method is the evaluation
of the joint $y$- and $S$- update step, which requires us to solve
a minimization problem that is not too much easier than the original
dual problem (\ref{eq:SDD}). ADMM overcomes this difficulty by adopting
an alternating-directions approach, updating $y$ while holding $S$
fixed, then updating $S$ using the new value of $y$ computed, as
in 
\begin{align*}
y^{k+1} & =\arg\min_{y}\Lag_{t}(X^{k},y,S^{k})\\
S^{k+1} & =\arg\min_{S\succeq0}\Lag_{t}(X^{k},y^{k+1},S)\\
X^{k+1} & =X^{k}+t\left(\sum_{i=1}^{m}y_{i}^{k+1}A_{i}+S^{k+1}-C\right).
\end{align*}
Here, the $y$- update is the unconstrained minimization of a quadratic
objective, and has closed-form solution 
\[
y^{k+1}=(\A^{T}\A)^{-1}\left[\frac{1}{t}\left(b-\A^{T}\vector X^{k}\right)+\A^{T}\vector(C-S^{k})\right],
\]
where $\A=[\vector A_{1},\ldots,\vector A_{m}]$. Similarly, the $S$-update
is the projection of a specific matrix matrix onto the positive-semidefinite
cone 
\[
S^{k+1}=[D]_{+}\text{ where }D=C-\sum_{i=1}^{m}y_{i}^{k+1}A_{i}-\frac{1}{t}X^{k},
\]
and also has a closed-form solution in terms of the eigenvalue decomposition
\begin{align*}
D & =\sum_{i=1}^{n}d_{i}vv^{T}, & [D]_{+} & =\sum_{i=1}^{n}\max\{d_{i},0\}vv^{T}.
\end{align*}
The iterates converge towards to the solutions of (\ref{eq:SDP})
and (\ref{eq:SDD}) for all fixed $t>0$~\cite[Theorem 2]{wen2010alternating},
although the sequence now typically converge much more slowly. In
practice, a heuristic based on balancing the primal and dual residuals
seems to work very well; see~\cite[Section 3.2]{wen2010alternating}
or \cite[Section 3.4.1]{boyd2011distributed} for its implementation
details.

\subsubsection{Complexity}

Unfortunately, it is difficult to bound the convergence rate of ADMM.
It has been shown that the sequence converges with sublinear objective
error $O(1/k)$ in an ergodic sense~\cite{he20121}, so in the worst
case, the method converges to $L$ accurate digits in $O(\exp(L))$
iterations. (This exponential factor precludes ADMM from being a polynomial-time
algorithm.) In practice, ADMM often performs much better than the
worst-case, converging to $L$ accurate digits in just $O(L)$ iterations
for a large array of SDP test problems~\cite{wen2010alternating}.

The per-iteration cost of ADMM can be dominated by the $y$-update,
due to the solution of the following system 
\begin{equation}
(\A^{T}\A)y=r\label{eq:ata}
\end{equation}
with a different right-hand side at each iteration. The standard approach
is to precompute the Cholesky factor, and to solve each instance (\ref{eq:ata})
by solving two triangular systems via forward- and back-substitution.
When the matrix $\A$ is fully-dense, the worst-case complexity of
solving (\ref{eq:ata}) is $O(n^{6})$ time and $O(n^{4})$ memory.

In practice, the matrix $\A$ is usually large-and-sparse, and complexity
of (\ref{eq:ata}) can be dramatically reduced by exploiting sparsity.
Indeed, the problem of solving a large-and-sparse symetric positive
definite system with multiple right-hand sides is classical in numerical
linear algebra, and the associated literature is replete. Efficiency
can be substantially improved by reordering the columns of $\A$ using
a fill-minimizing ordering like minimum degree and nested dissection~\cite{george1981computer},
and by using an incomplete factorization as the preconditioner within
an iterative solution algorithm like conjugate gradients~\cite{saad2003iterative}.
The cost of solving practical instances of (\ref{eq:ata}) can be
as low as $O(m)$.

Assuming that the cost of solving (\ref{eq:ata}) can be made negligible
by exploiting sparsity in $\A$, the per-iteration cost is then dominated
by the eigenvalue decomposition required for the $S$-update. Performing
this step using dense linear algebra requires $\Theta(n^{3})$ time
and $\Theta(n^{2})$ memory. For larger problems where $X^{\star}$
is known to be low-rank, it may be possible to use low-rank linear
algebra and an iterative eigendecomposition like Lanczos iterations
to push the complexity figure down to as low as $O(n)$ per-iteration.

\subsubsection{Bibliography}

ADMM was originally proposed in the mid-1970s by Glowinski and Marrocco~\cite{glowinski1975approximation}
and Gabay and Mercier~\cite{gabay1976dual}, and was studied extensively
in the context of maximal monotone operators (see~\cite[Section 3.5]{boyd2011distributed}
for a summary of the historical developments). The algorithm experienced
a revival in the past decade, in a large part due to the publication
of a popular and influential survey by Boyd~et~al.~\cite{boyd2011distributed}
for applications in distributed optimization and statistical learning.
The algorithm described in this subsection was first proposed by Wen,
Goldfarb and Yin~\cite{wen2010alternating}, and is one of two popular
variations of ADMM specifically designed for the solution of large-scale
SDPs, alongside the algorithm of O'Donoghue~et.~al~\cite{o2016conic}.

\subsection{\label{subsec:mIPM}Modified interior-point method for low-rank SDPs}

A fundamental issue with standard off-the-shelf interior-point solvers
is their inability to exploit problem structure to substantially reduce
complexity. In this subsection, we describe a modification to the
standard interior-point method that makes it substantially more efficient
for large-and-sparse low-rank SDPs, for which the number of nonzeros
in the data $A_{1},\ldots,A_{m}$ is small, and $\theta\triangleq\rank\{X^{\opt}\}$
is known \emph{a priori} to be very small relative to the dimensions
of the problem, i.e. $\theta\ll n$. As we have previous reviewed
in the previous two sections, such problems widely appear by applying
convexification to problems in graph theory~\cite{javadsojoudi14},
approximation theory~\cite{candes2010power,lasserre-2010}, control
theory~\cite{valmorbida2016stability,lasserre-2010,zhang2016thesis},
and power systems~\cite{LL2012_1}. They are also the fundamental
building blocks for global optimization techniques based upon polynomial
sum-of-squares~\cite{parrilo-2003} and the generalized problem of
moments~\cite{lasserre-2010}.

To describe this modification, let us recall that modern primal-dual
interior-point methods almost always converge in 30-50 iterations,
and that their per-iteration cost is dominated by the solution of
the Hessian equation 
\begin{equation}
\underbrace{[\A^{T}(W\otimes W)\A]}_{\H}y=r,\label{eq:hesseqn1_2}
\end{equation}
in which $\A=[\vector A_{1},\ldots,\vector A_{m}]$, and $W$ is the
positive definite scaling matrix. An important feature of interior-point
methods for SDPs is that the matrix $W$ is fully-dense, and this
makes $\H$ fully-dense, despite any apparent sparsity in the data
matrix $\A$. The standard approach of dense Cholesky factorization
takes approximately the same amount of time and memory for sparse,
low-rank problems as it does for dense, high-rank problems.

Alternatively, the Hessian equation may be solved iteratively using
the preconditioned conjugate gradients (PCG) algorithm. We defer to
standard texts~\cite{barrett1994templates} for the implementation
details of PCG, and only note that at each iteration, the method requires
a single matrix-vector product with the governing matrix $\H$, and
a single \emph{solve}\footnote{i.e. matrix-vector product with the inverse $\tilde{\H}^{-1}$.}
with the preconditioner $\tilde{\H}$. The key ingredient is a good
preconditioner: a matrix $\tilde{\H}$ that is similar to $\H$ in
a spectral sense, but is otherwise much cheaper to invert. The following
iteration bound is standard; a proof can be found in standard references
like~\cite{saad2003iterative} and~\cite{greenbaum1997iterative}. 
\begin{prop}
\label{prop:CG}Consider using preconditioned conjugate gradients
to solve $\H y=r$, with $\tilde{\H}$ as preconditioner. Define $y^{\star}=\H^{-1}r$
as the exact solution, and $\kappa=\lambda_{\max}(\tilde{\H}^{-1}\H)/\lambda_{\min}(\tilde{\H}^{-1}\H)$
as the joint condition number. Then at most 
\[
i\le\left\lceil \frac{\sqrt{\kappa}}{2}\log\left(\frac{2\sqrt{\kappa}}{\epsilon}\right)\right\rceil \text{ PCG iterations}
\]
are required to compute an $\epsilon$-accurate iterate $y^{i}$ satisfying
$\|y^{i}-y^{\star}\|\le\epsilon\|y^{\star}\|$. 
\end{prop}
A preconditioner that guarantees joint condition number $\kappa=O(1)$
was described in~\cite{zhang2017modified}. The preconditioner based
on the insight that the scaling matrix $W$ can be decomposed into
two components, 
\begin{equation}
W=W_{0}+UU^{T},\label{eq:WUUt}
\end{equation}
in which the rank of $U$ is at most $\theta=\rank\{X^{\opt}\}$,
and $W_{0}$ is well-conditioned, meaning that all of its eigenvalues
are roughly the same value. Substituting (\ref{eq:WUUt}) into (\ref{eq:hesseqn1_2})
reveals the same decomposition for the Hessian matrix, 
\begin{equation}
\H=\A^{T}(W_{0}\otimes W_{0})\A+\underbrace{\A^{T}(U\otimes Z)(U\otimes Z)^{T}\A}_{\U\U^{T}}\label{eq:hessExp2}
\end{equation}
where $Z$ is any matrix (not necessarily unique) satisfying $ZZ^{T}=2E+UU^{T}$.
Note that the rank of $\U$ is at most $n\theta$.

To develop a preconditioner based on this insight, we make the approximation
$W_{0}\approx\tau I$ in (\ref{eq:hessExp2}), to yield 
\begin{equation}
\tilde{\H}=\tau^{2}\A^{T}\A+\U\U^{T}.\label{eq:AtIUUtAt}
\end{equation}
This dense matrix is a low-rank perturbation of the sparse matrix
$\A^{T}\A$, and can be inverted using the Sherman–Morrison–Woodbury
formula 
\begin{equation}
\tilde{\H}^{-1}=(\tau^{2}\A^{T}\A)^{-1}(I-\U\Sch^{-1}\U^{T}(\A^{T}\A)^{-1}),\label{eq:SMW2}
\end{equation}
in which the Schur complement $\Sch=\tau^{2}I+\U^{T}(\A^{T}\A)^{-1}\U$
can be precomputed. 
\begin{lem}[{{\cite[Lemma 7]{zhang2017modified}}}]
\label{lem:Htilde}Let $\tilde{\H}$ be defined in (\ref{eq:AtIUUtAt}),
and choose $\tau$ to satisfy $\lambda_{\min}(E)\le\tau\le\lambda_{\max}(E)$.
Then, the joint condition number $\kappa=\lambda_{\max}(\tilde{\H}^{-1}\H)/\lambda_{\min}(\tilde{\H}^{-1}\H)$
is an absolute constant. 
\end{lem}
Combining Proposition~\ref{prop:CG} and Lemma~\ref{lem:Htilde},
we find that PCG with $\tilde{\H}$ as preconditioner will solve the
Hessian equation to $L$ digits of accuracy in at most $O(L)$ iterations.

\subsubsection{Complexity}

All interior-point methods converge to $L$ accurate digits in between
$O(\sqrt{n}L)$ and $O(nL)$ Newton iterations, and practical implementations
almost always occurs with 30-50 Newton iterations. Performing each
Newton iteration using the PCG procedure described above, this translates
into a formal complexity bound of $O(\sqrt{n}L^{2})$ to $O(nL^{2})$
PCG iterations, and a practical value of between 500-1500 PCG iterations.

The per-iteration cost of PCG can be dominated by the cost of applying
the Sherman–Morrison–Woodbury formula to invert the preconditioner
$\tilde{\H}$, due to the need of repeatedly making matrix-vector
products with $(\A^{T}\A)^{-1}$. Like discussed in Section~\ref{subsec:ADMM}
for ADMM, the matrix $\A$ is usually large-and-sparse in practical
applications, and standard techniques from numerical linear algebra
can often substantially reduce the cost of this operation. Assuming
that this matrix-vector product can be performed in $O(m)$ time and
memory, \cite{zhang2017modified} showed that the cost of inverting
$\tilde{\H}$ is $O(\theta^{2}n^{2})$ time and memory, after an initial
factorization step requiring $O(\theta^{3}n^{3})$ time, where $\theta=\rank\{X^{\opt}\}$.

Assming that the cost of applying the preconditioner $\tilde{\H}$
is negligible, the per-iteration cost of PCG becomes dominated by
the matrix-vector with $\H$, which is $\Theta(n^{3})$ time and $\Theta(n^{2})$
memory. Assuming that $\theta$ and $L$ are both significantly smaller
than $n$, the formal complexity of the algorithm is $\Theta(n^{3.5})$
time and $\Theta(n^{2})$ memory, and the practical complexity is
closer to $\Theta(n^{3})$ time.

\subsubsection{Bibliography}

The idea of using conjugate gradients (CG) to solve the Hessian equation
dates back to the original Karmarkar interior-point method~\cite{karmarkar1984new},
and was widely used in early interior-point codes for SDPs~\cite{vandenberghe1995primal,vandenberghe1996semidefinite}.
However, subsequent numerical experience~\cite{fujisawa2000numerical,benson2001dsdp3}
found the approach to be ineffective: the Hessian matrix $\H$ becomes
ill-conditioned as the interior-point iterate approaches the solution,
and CG requires more and more iterations to converge. Toh and Kojima~\cite{toh2002solving}
were the first to develop highly effective \emph{spectral} preconditioners,
based on the same decomposition of the scaling matrix $W$ as above.
However, its use required almost as much time and memory as a single
iteration of the regular interior-point method. The modified interior-point
method of~\cite{zhang2017modified}, which we had described in this
subsection, makes the same idea efficient by utilizing the Sherman–Morrison–Woodbury
formula.

\subsection{\label{subsec:other}Other specialized algorithms}

This tutorial has given an overview of the interior-point method as
a general-purpose algorithm for SDPs, and described two specialized
structure-exploiting algorithms for large-scale SDPs in detail. Numerous
other structure-exploiting algorithms also exist. In general, it is
convenient to categorize them into three distinct groups:

The first group comprises first-order methods, like ADMM in Section~\ref{subsec:ADMM},
but also smooth gradient methods~\cite{nesterov2007smoothing}, conjugate
gradients~\cite{toh2002solving,zhao2010newton,zhang2016thesis},
augmented Lagrangian methods~\cite{kovcvara2003pennon}, applied
either to (\ref{eq:SDP}) directly, or to the Hessian equation associated
with an interior-point solution of (\ref{eq:SDP}). All of these algorithms
have inexpensive per-iteration costs but a sublinear worst-case convergence
rate, computing an $\epsilon$-accurate solution in $O(1/\epsilon)$
time. They are most commonly used to solve very large-scale SDPs to
modest accuracy, though in fortunate cases, they can also converge
to high accuracy.

The second group comprises second-order methods that use sparsity
in the data to decompose the size-$n$ conic constraint $X\succeq0$
into many smaller conic constraints over submatrices of $X$. In particular,
when the matrices $C,A_{1},\ldots,A_{m}$ share a common sparsity
structure with a chordal graph with bounded treewidth $\tau$, a technique
known as \emph{chordal decomposition} or \emph{chordal conversion}
can be used to reformulate (\ref{eq:SDP})-(\ref{eq:SDD}) into a
problem containing only size-$(\tau+1)$ semidefinite constraints~\cite{fukuda2001exploiting};
see also~\cite{vandenberghe2015chordal}. While the technique is
only applicable to chordal SDPs with bounded treewidths, it is able
to reduce the cost of a size-$n$ SDP all the way down to the cost
of a size-$n$ linear program, sometimes as low as $O(\tau^{3}n)$.
Indeed, chordal sparsity can be guaranteed in many important applications~\cite{vandenberghe2015chordal,Madani_SIAM_2017},
and software exist to automate the chordal reformulation~\cite{kim2011exploiting}.

The third group comprises formulating low-rank SDPs as nonconvex optimization
problems, based on the outer product factorization $X=RR^{T}$. The
number of decision variables is dramatically reduced from $\sim n^{2}$
to $n$~\cite{burer2003nonlinear,journee2010low}, though the problem
being solved is no longer convex, so only local convergence can be
guaranteed. Nevertheless, time and memory requirements are substantially
reduced, and these methods have been used to solve very large-scale
low-rank SDPs to excellent precision; see the computation results
in~\cite{burer2003nonlinear,journee2010low}.

\section{Conclusion}

\label{sec:5}

Optimization lies at the core of classical control theory, as well
as up-and-coming fields of statistical and machine learning. This
tutorial paper provides an overview of conic optimization, and its
application to the design, analysis, control and operation of real-world
systems. In particular, we give concrete case studies on machine learning,
power systems, and state estimation, as well as the abstract problems
of rank minimization and quadratic optimization. We show that a wide
range of nonconvex problems can be converted in a principled manner
into a hierarchy of convex problems, using a range techniques collectively
known as convexification. Finally, we develop numerical algorithms
to solve these convex problems in a highly efficient manner, by exploiting
problem structure like sparsity and low solution rank.

\bibliographystyle{IEEEtran}
\bibliography{combined}

\end{document}